\documentclass[11pt]{article}

\usepackage[T1]{fontenc}
\usepackage{lmodern}
\usepackage[a4paper,margin=30mm]{geometry}
\usepackage{amsmath,amssymb,amsthm}
\usepackage[numbers,sort&compress]{natbib}
\usepackage[section]{placeins}
\usepackage{flafter}
\usepackage{microtype}
\usepackage[hidelinks]{hyperref}

\hypersetup{
  pdftitle={Finite tests and low-rank excluded minors for matroid path-width two},
  pdfauthor={Antonios Kalampakas},
  pdfsubject={Matroid path-width and excluded minors},
  pdfkeywords={matroid path-width, excluded minor, binary matroid,
    two-sum, parallel extension, computational classification}
}

\newtheorem{theorem}{Theorem}[section]
\newtheorem{lemma}[theorem]{Lemma}
\newtheorem{proposition}[theorem]{Proposition}
\newtheorem{corollary}[theorem]{Corollary}
\theoremstyle{definition}
\newtheorem{definition}[theorem]{Definition}

\theoremstyle{remark}

\DeclareMathOperator{\pw}{pw}

\DeclareMathOperator{\si}{si}

\newcommand{\F}{\mathbb{F}_2}

\title{Finite tests and low-rank excluded minors\\
for matroid path-width two}

\author{Antonios Kalampakas\\[0.5em]
\small College of Engineering and Technology,\\
\small American University of the Middle East,\\
\small Egaila 54200, Kuwait\\[0.3em]
\small\href{mailto:antonios.kalampakas@aum.edu.kw}{antonios.kalampakas@aum.edu.kw}\\
\small ORCID: \href{https://orcid.org/0000-0001-8821-6102}{0000-0001-8821-6102}}

\date{}

\begin{document}
\maketitle

\begin{abstract}
The path-width of a matroid is the minimum, over all orderings of its ground set, of the maximum value of its connectivity function on an initial segment. For any two connected rooted matroids whose roots are neither loops nor coloops, we form the 2-sum of each with a third such rooted matroid by identifying their roots and deleting the common element. The two original rooted matroids are equivalent when these sums agree on whether their path-width is at most two for every choice of the third rooted matroid. We prove that the two rooted matroids are equivalent if and only if this agreement holds for thirteen fixed graphic or cographic rooted matroids. They distinguish exactly twenty equivalence classes, and, moreover, there exists no family of twelve fixed connected rooted matroids that can guarantee this. Using these classes, we show that the tree in the canonical tree decomposition of an excluded minor for path-width at most two has diameter at most ten edges. For a connected matroid given by a tree decomposition using 2-sums, we give an exact recurrence to decide whether its path-width is at most two and to construct an ordering of width at most two when possible. Moreover, we derive a formula for the path-width of parallel extensions that gives complete binary excluded-minor lists in ranks three, four and five, containing two, four and eleven isomorphism classes, respectively.
\end{abstract}


\emph{Secondary:} 05C30, 68R05.

\section{Introduction}\label{sec:intro}
\begingroup
\makeatletter
\@clubpenalty=10000
\makeatother
\clubpenalty=10000
\widowpenalty=10000

For a matroid $M$ with ground set $E$ and rank function $r$, the connectivity function is $\lambda_M(X)=r(X)+r(E\setminus X)-r(E)$ for $X\subseteq E$. The width of an ordering $e_1,\ldots,e_n$ of $E$ is the largest value of $\lambda_M(\{e_1,\ldots,e_i\})$ for $1\leq i\leq n$, and the path-width of $M$ is the minimum of this width over all orderings. This parameter was introduced by Geelen, Gerards and Whittle~\cite{GGW06}. For a matroid represented over a finite field, Kashyap identifies its path-width with the minimum trellis state complexity of the linear code specified by the representation~\cite{Kas08}.

We study how path-width at most two behaves under $2$-sums and what this implies for excluded minors. A \emph{rooted matroid} $(N,p)$ is a matroid $N$ with a distinguished element $p$, called its \emph{root}. We require the root of every rooted matroid considered below to be neither a loop nor a coloop. For a $2$-sum $N\oplus_2 C$, the roots of the summands are identified as a common element $p$, and the ground set of the sum is $(E(N)\cup E(C))\setminus\{p\}$~\cite{Oxl11}. To compare connected rooted matroids $(N,p)$ and $(N',p)$, we form their $2$-sums with the same connected rooted matroid $(C,p)$, taking the remaining ground sets disjoint. Our first question is when $(N,p)$ and $(N',p)$ are equivalent, meaning that the two sums either both have path-width at most two or both have greater path-width for every connected rooted matroid $(C,p)$.

The class of matroids of path-width at most two is minor-closed, so we also study its excluded minors, i.e., matroids of greater path-width whose proper minors all have path-width at most two. For matroids representable over a fixed finite field, there are only finitely many such excluded minors~\cite{GGW02,KKKO23}. Kashyap gave a partial list at width two, while Koutsonas, Thilikos and Yamazaki classified the outerplanar obstructions~\cite{Kas08,KTY14}. The complete binary list at width two remains unknown. We derive restrictions on the $2$-sum decompositions of excluded minors and determine the complete binary lists through rank five, without imposing a bound on the number of elements.

For the replacement question, we prove that thirteen fixed graphic or cographic rooted matroids suffice, each with at most fourteen elements. For each of these thirteen matroids, form its $2$-sum with $N$ and with $N'$ by identifying their roots. The rooted matroids $N$ and $N'$ are equivalent if and only if, for each of the thirteen choices, the two resulting sums agree on whether their path-width is at most two. There are exactly twenty equivalence classes, each with a graphic or cographic representative, and no family of twelve fixed connected rooted matroids can distinguish all of them. Thus thirteen is optimal. Both the equivalence criterion and the optimality statement apply to all connected rooted matroids, without assuming representability.

To prove the criterion, we specialise the trajectory and replacement methods of Jeong, Kim and Oum and of Kant\'e, Kim, Kwon and Oum~\cite{JKO17,KKKO23}. Elements from different summands can be interleaved in an ordering of a $2$-sum, so the separate path-widths of the summands do not retain all the information needed for this comparison. At each cut of an ordering of $E(N)\setminus\{p\}$, we retain the sum of the ranks of its two parts minus $r(N)$, together with whether $p$ lies in the closure of each part. We compress these sequences while preserving whether they can be interleaved with another sequence to give an ordering of width at most two. The compressed sequences lie in an explicit set of 114 possibilities, each of length at most nine. We construct a table specifying which pairs of these sequences can be interleaved to give an ordering of width at most two, and this table gives the twenty equivalence classes and the thirteen fixed matroids.

For a connected matroid given by a tree decomposition using $2$-sums, we give an exact recurrence to decide whether its path-width is at most two and to construct an ordering of width at most two when possible. The vertices of the decomposition tree are labelled by connected matroids, and each edge specifies a $2$-sum at a distinguished element common to the matroids labelling its ends~\cite{Oxl11}. The distinguished elements used for different edges incident with the same vertex are distinct. The recurrence combines families of compressed sequences for the matroids labelling the vertices, which may include circuits and cocircuits. Since a step in a compressed sequence may span several original elements, we use a fixed containing set of 208 sequences during the intermediate calculations. We check the decisions and every cut of the constructed orderings independently on the resulting matroids. There are already algorithms for constructing orderings of small width for matroids represented over a fixed finite field~\cite{Hli18,JKO17}. Here we give explicit rules for tree decompositions using $2$-sums, without claiming polynomial running time when the sizes of the matroids labelling the vertices or the degrees of the vertices in the decomposition tree are unrestricted.

A simpler criterion holds when one summand $N$ is $3$-connected and has at least four elements. If the other summand is connected and is neither a circuit nor a cocircuit, every $2$-sum of path-width at most two admits an ordering of width at most two in which $E(N)\setminus\{p\}$ is consecutive. The replacement question for $N$ therefore reduces to five conditions, concerning an ordering with $p$ at an end, deletion of $p$, contraction of $p$, a series extension at $p$ and a parallel extension at $p$. We obtain at most eight equivalence classes for these $3$-connected summands. Our endpoint arguments use results of Hall, Oxley and Semple~\cite{HOS07} and of Beavers and Oxley~\cite{BO08}, whose width convention is one unit higher than ours.

Now consider a connected matroid given by a tree decomposition using $2$-sums in which every vertex is labelled by a $3$-connected matroid on at least four elements. As above, the $2$-sums use a distinct distinguished element for each edge incident with a vertex. The resulting matroid has path-width at most two if and only if the tree is a path and each matroid labelling a vertex has an ordering of width at most two that places its distinguished elements at the ends. At an internal vertex of the path, the same ordering must begin with one distinguished element and end with the other, so separate orderings placing each element at an end do not suffice. Our examples show that, for a path with three vertices, the $2$-sums of the matroids labelling its two pairs of adjacent vertices can both have path-width two even though the sum of all three matroids has path-width three. Circuits and cocircuits can label vertices in decomposition trees that are not paths even when the resulting matroid has path-width at most two, so the general recurrence retains the compressed cut sequences.

Using the criterion based on the twenty equivalence classes, we show that the tree in the canonical tree decomposition of every excluded minor for path-width at most two has diameter at most ten edges. Every circuit or cocircuit labelling a vertex of this decomposition has at most eleven elements, counting the distinguished elements used for the $2$-sums, and both bounds apply to arbitrary matroids. The proof uses minor operations to shorten a path and compare the results of the $2$-sums with the thirteen fixed matroids on its two sides. At each cut, the table of pairs of sequences described above assigns an integer from one to ten which increases strictly along the path, giving the diameter bound.

For the binary excluded minors, the complete lists contain two isomorphism classes of rank three, four of rank four and eleven of rank five. The maximum numbers of elements at ranks four and five are twelve and eighteen, respectively. The $3$-sun consists of a triangle with a distinct additional vertex joined to the endpoints of each triangle edge. Doubling every edge and taking the cycle matroid gives the unique rank-five class with eighteen elements. The doubled graph is in the outerplanar obstruction family classified in~\cite[Theorem~1]{KTY14}, while our complete rank-five list shows that no other class has eighteen elements.

To obtain these classifications, we calculate the path-width after parallel extensions by an exact formula, proved using the reordering argument in~\cite[proof of Lemma~3.8]{Kas08}. With Tutte's graphic characterization, the formula reduces the graphic searches to doubling subsets of edges of ten simple graphs for rank four and fifty-six for rank five. We verify that every required simple graph is included, repeat the path-width calculations and verify which candidates are excluded minors. The lists are therefore complete at these ranks, and the bounds of twelve and eighteen elements follow from the classification, without imposing a size limit on the candidates.

\par
\endgroup

\section{Preliminaries}\label{sec:prelim}

All graphs and matroids are finite. Graphs are simple unless described as
multigraphs, which have no loops unless stated. A $2$-connected multigraph
is connected, has at least two edges, and has no cut vertex, so two
vertices joined by parallel edges are allowed. We follow Oxley~\cite{Oxl11} for matroid terminology. In particular,
$E(M)$ and $r_M$ denote the ground set and rank function of a matroid $M$,
and $r(M)=r_M(E(M))$. We write $M^*$ for its dual, $M\setminus e$
for deletion and $M/e$ for contraction. The closure of $X$ is
$\operatorname{cl}_M(X)=\{e\in E(M):r_M(X\cup\{e\})=r_M(X)\}$.
An indicator $\mathbf1_P$, also written $[P]$, is one when the
statement $P$ holds and zero otherwise.

The results about $2$-sums apply to arbitrary matroids. The excluded-minor
classifications concern binary matroids, i.e., matroids represented by
matrices over $\F$, the field with two elements. All matrix ranks in those
classifications are taken over $\F$.

We use the connectivity function \[
\lambda_M(X)=r_M(X)+r_M(E(M)\setminus X)-r(M),
\] as in \cite{Kas08,KTY14,KKKO23}, whereas the conventions of \cite{GGW06,HOS07} are one unit higher.

\begin{definition}\label{def:pw}
For an ordering $e_1,\ldots,e_n$ of $E(M)$, a cut separates an initial segment
$\{e_1,\ldots,e_i\}$ from the remaining elements. The width of the ordering is
\[
\max_{1\leq i\leq n}\lambda_M(\{e_1,\ldots,e_i\}).
\]
The \emph{path-width} $\pw(M)$ of $M$ is the minimum of these widths over all orderings of $E(M)$, with value zero for the empty matroid.
\end{definition}

Deleting or contracting elements cannot increase path-width. Moreover,
$\lambda_M=\lambda_{M^*}$, so $\pw(M)=\pw(M^*)$. We also use the bound
$\lambda_M(X\cup\{e\})\leq\lambda_M(X)+1$, which follows because
adding one element increases rank by at most one.

The direct-sum formula
$\pw(M_1\oplus M_2)=\max\{\pw(M_1),\pw(M_2)\}$ follows by
concatenating orderings and applying minor monotonicity.
We take width bounds $k$ to be nonnegative integers and call a matroid
an \emph{excluded minor} for path-width at most $k$ when its path-width is greater
than $k$ and every proper minor has path-width at most $k$.

\begin{proposition}\label{prop:exact-value}
Every excluded minor for path-width at most $k$ is connected and has
path-width exactly $k+1$.
\end{proposition}

\begin{proof}
If the matroid were disconnected, the direct-sum formula would give
a proper component of path-width greater than $k$, contradicting minimality.
Now delete an element $e$ and choose an ordering of $M\setminus e$
of width at most $k$. Appending $e$ increases the width by at most one,
since $\lambda_M(X)\leq\lambda_{M\setminus e}(X)+1$ for each old
initial segment $X$. The excluded minor therefore has path-width at most $k+1$,
and its path-width is greater than $k$ by definition.
\end{proof}

\section{Parallel extensions}\label{sec:structure}

All results in this section hold without representability assumptions.
For a matroid $N$, write $N_\mu$ for the matroid obtained by replacing
each nonloop $e$ by a parallel class of size $\mu(e)\geq1$, while leaving
the loops unchanged. For $D\subseteq E(N)$, let $N[D]$ add one parallel copy
of every element of $D$, with copies of loops again taken as loops.

To calculate the path-width after adding parallel copies, we first consider
an ordering in which each parallel class is consecutive. A cut between classes has the same connectivity as the corresponding cut in $N$, whereas a cut inside the class of $e$ has a copy of $e$ on each side and therefore requires a separate rank expression. For an ordering $\pi$ of $E(N)$,
let $A_\pi(e)$ and $B_\pi(e)$ be the sets of elements before and after
$e$, respectively, and put
\[
 \chi_\pi(e)=r_N(A_\pi(e)\cup\{e\})
             +r_N(B_\pi(e)\cup\{e\})-r(N).
\]
Write $w_N(\pi)$ for the width of $\pi$. For $e\notin A$, define
\[
 \chi(A,e)=r_N(A\cup\{e\})+r_N(E(N)\setminus A)-r(N),
\]
so $\chi_\pi(e)=\chi(A_\pi(e),e)$.

\begin{theorem}[Exact parallel-extension formula]\label{thm:exact-doubling}
For every matroid $N$ and $D\subseteq E(N)$,
\[
 \pw(N[D])=
 \min_\pi\max\bigl(\{w_N(\pi)\}\cup
                         \{\chi_\pi(e):e\in D\}\bigr).
\]
The same formula holds for arbitrary parallel multiplicities, with $D$
the set of elements having multiplicity at least two.
\end{theorem}

Using the reordering argument in Kashyap's proof of
Lemma~3.8~\cite{Kas08}, we show that the parallel classes can be
made consecutive without increasing width.

\begin{proof}
In any ordering of the parallel extension, let $e$ be the first member of a designated
parallel class and move a later member $f$ immediately after $e$. Every changed initial segment has the form $P\cup\{f\}$ for an old initial segment $P$ containing $e$. Since $e$ and $f$ are parallel, this initial segment has rank $r(P)$, and removing $f$ from its complement cannot increase the complement's rank, so the move does not increase width. Repeating the operation for each designated class, while preserving the classes already grouped, therefore gives an optimal ordering in which all designated classes are consecutive.

The induced ordering $\pi$ of $E(N)$ has the ordinary connectivity values
between classes. Inside the class of $e$, both sides contain a representative
of $e$, so the connectivity is exactly $\chi_\pi(e)$. This proves the
formula and the multiplicity statement. Repeated loops have no effect.
\end{proof}

The formula also gives a dynamic program on the subsets of $E(N)$.
For a partial ordering of $S$, its cost is the maximum of the connectivity values of its initial segments and the values of $\chi(A,e)$ encountered when $e\in D$ is appended after the initial segment $A$, with all ranks and complements evaluated in the full matroid $N$. Let $F(S)$ be the minimum of
this cost over the orderings of $S$. Put $b_D(A,e)=\chi(A,e)$ when
$e\in D$, and $b_D(A,e)=0$ otherwise. Then $F(\varnothing)=0$ and
\[
 F(S)=\min_{e\in S}\max\{
 F(S\setminus\{e\}),\lambda_N(S),b_D(S\setminus\{e\},e)\}.
\]
Indeed, choosing $e$ as the last element of $S$ adds the cut at $S$
and, when $e\in D$, the cut inside its parallel class. Minimizing over
that choice gives the recurrence, and Theorem~\ref{thm:exact-doubling}
gives $F(E(N))=\pw(N[D])$. If $n=|E(N)|$ and all subset ranks are supplied, the recurrence has $O(n2^n)$ transitions, independently of how many parallel copies are added to $N$. For represented matroids, $\chi_\pi(e)$ is the trellis
branch-complexity at that position \cite[Eq.~(4)]{KT12}.

\begin{lemma}[Consequences for parallel extensions]
\label{thm:doubling}\label{thm:parclass}\label{lem:loops-coloops}
For every parallel extension $N_\mu$ of $N$,
\[
 \pw(N)\leq\pw(N_\mu)\leq\pw(N)+1.
\]
The dual bound holds for series extensions. Every excluded minor for
path-width at most $k\geq1$ has no loops or coloops, and every parallel
and series class has at most two elements.
\end{lemma}

\begin{proof}
The lower bound is minor monotonicity, and
$\chi_\pi(e)\leq\lambda_N(A_\pi(e))+1$ gives the upper bound by
Theorem~\ref{thm:exact-doubling}. The formula depends only on which
multiplicities exceed one, so deleting a third member of a parallel
class leaves path-width unchanged. Inserting or deleting a loop also
leaves all connectivity values unchanged. These observations rule out
loops and parallel classes of size at least three in an excluded minor.
Duality gives the statements about coloops and series classes.
\end{proof}

The simplification $\si(M)$ deletes loops and all but one element of
each parallel class. Since Lemma~\ref{thm:parclass} bounds every parallel class of an excluded minor by two, every loopless excluded minor $M$ for path-width at most two has the form $N[D]$, where $N=\si(M)$ is simple and $D\subseteq E(N)$.

\section{Two-sums of three-connected matroids and endpoint orderings}\label{sec:two-sums}

An ordering of width at most two of a $3$-connected matroid need not permit an arbitrary
element at an end. The endpoint structure of such orderings is treated in
\cite[Theorems~1.2--1.4 and~2.3]{HOS07}. The following proposition links
endpoint conditions in two $3$-connected summands to the path-width of
their $2$-sum. The rooted summands meet in their common root $p$, which is neither a loop nor a coloop of either summand. Their other elements are disjoint, and $p$ is removed when the $2$-sum is formed~\cite{Oxl11}. Thus, writing $M=N_1\oplus_2N_2$ and $A_i=E(N_i)\setminus\{p\}$, we have $E(M)=A_1\mathbin{\dot\cup}A_2$. An ordering lists $A_i$ consecutively
when its elements form an interval of that ordering.

For a matroid $N$ and $p\in E(N)$, let $\pw_p(N)$ denote the minimum width of an ordering with $p$ at an end, where either choice of end gives the same parameter because reversal preserves width.

\begin{proposition}[Two-sums of $3$-connected matroids]
\label{prop:rooted-two-sum}
Let $(N_1,p)$ and $(N_2,p)$ be rooted $3$-connected matroids, each with at least four elements, and let $M=N_1\oplus_2N_2$.
Then
\[
  \pw(M)\leq2
  \quad\Longleftrightarrow\quad
  \pw_p(N_1)\leq2\ \text{ and }\ \pw_p(N_2)\leq2.
\]
Every ordering of $M$ of width at most two lists the two sets
$E(N_1)\setminus\{p\}$ and $E(N_2)\setminus\{p\}$ consecutively.
\end{proposition}

\begin{proof}
Put $A_i=E(N_i)\setminus\{p\}$. For $X_i\subseteq A_i$, write
$Y_i=A_i\setminus X_i$ and define
\[
 a_i=\mathbf1_{p\in\operatorname{cl}_{N_i}(X_i)},\qquad
 b_i=\mathbf1_{p\in\operatorname{cl}_{N_i}(Y_i)},\qquad
 d_i=r_{N_i}(X_i)+r_{N_i}(Y_i)-r(N_i).
\]
Here $a_i$ indicates whether $p$ is in the closure of the listed part $X_i$, and $b_i$ indicates whether it is in the closure of the unlisted part $Y_i$. These two indicators determine whether the common element
reduces the rank when the corresponding parts are joined.
The rank formula for a $2$-sum gives
\begin{equation}\label{eq:two-sum-cut}
 \lambda_M(X_1\cup X_2)=d_1+d_2+1-a_1a_2-b_1b_2.
\end{equation}
If $X_i$ is a singleton, simplicity and cosimplicity of $N_i$ give $(d_i,a_i,b_i)=(1,0,1)$, whereas if $Y_i$ is a singleton the tuple is $(1,1,0)$. If both $X_i$ and $Y_i$ have at least two elements,
$3$-connectivity applied with $p$ on either side gives
\[
 d_i\geq1+\max\{a_i,b_i\}.
\]
Consequently, if both $A_1$ and $A_2$ contain listed and unlisted
elements, the cut value is at least three, except possibly when both
listed parts are singletons or both unlisted parts are singletons.
For the case in which all four parts have at least two elements, use
$a_1a_2+b_1b_2\leq\max\{a_1,b_1\}+\max\{a_2,b_2\}$ in
\eqref{eq:two-sum-cut}.

Consider the first switch between sides in an ordering of width at most two. If at
least two elements of the first side have been listed but that side is
incomplete, the switch has cut value at least three. If only one element has been listed, the switch leaves one listed element in each side, and the next element gives a cut of value at least three because $|A_1|,|A_2|\geq3$. Thus the first side must be completed before the
second begins.

For a concatenation of the two side orderings, the cuts during the first side have values $\lambda_{N_1}(X_1)$, with $p$ on the unlisted side, while those during the second have values $\lambda_{N_2}(X_2\cup\{p\})$. Such a concatenation therefore has minimum width $\max\{\pw_p(N_1),\pw_p(N_2)\}$, which, together with the preceding necessity, proves the equivalence.
\end{proof}

The size hypothesis is needed for the assertion about every ordering.
Two triangles have a four-element circuit as their $2$-sum. Every ordering of this circuit, including an interleaved ordering, has width one.
To apply the same endpoint condition to a tree decomposition using $2$-sums, let $T$ be a finite tree and assign to every vertex $v$ a $3$-connected matroid $N_v$ with at least four elements, using disjoint ground sets.
For every edge $uv$, choose one element in each of $N_u$ and $N_v$
as distinguished elements for the $2$-sum, with distinct choices for the edges incident with a fixed vertex. Identify the chosen pair and take the $2$-sum at that element. Write $M(T)$ for the resulting matroid and $A_v$ for the elements of $N_v$ not used in these $2$-sums. These are precisely the
elements of $N_v$ that remain in $M(T)$, so
\[
 E(M(T))=\mathop{\dot\bigcup}_{v\in V(T)} A_v.
\]
The order of these $2$-sums does not affect $M(T)$, and the sum over
any connected subtree is connected~\cite{Oxl11}. Deleting an edge $uv$ of $T$ gives two components. The unions of the sets $A_v$ over these components are the two parts of the partition of $E(M(T))$ displayed by that edge. For a component of $T-v$, we take the $2$-sums of the matroids labelling its vertices and retain the distinguished element for the edge joining it to $v$. This gives a rooted matroid. We use these rooted matroids when forming the $2$-sums with $N_v$.

For distinct elements $p,q$ of a matroid $N$, let $\pw_{p,q}(N)$ be
the minimum width of an ordering that starts with $p$ and ends with
$q$. Reversing the ordering exchanges $p$ and $q$.

\begin{lemma}[Consecutive parts of the partition displayed by an edge]\label{lem:tree-shore-blocks}
In every ordering of $M(T)$ of width at most two, the two parts of the partition displayed by each edge of $T$ are consecutive blocks, so one is an initial segment and the other is a final segment.
\end{lemma}

\begin{proof}
To reduce a connected rooted matroid $(B,p)$ while retaining its root $p$ and any chosen other element $e$, choose a circuit $C$ containing $p$ and $e$, which exists by connectedness. Contracting $C\setminus\{p,e\}$ and
deleting all remaining elements outside $\{p,e\}$ gives the parallel
pair $U_{1,2}$. Consequently, when $(B,p)$ is a summand, it can be replaced by any chosen remaining element $e$, which takes the place of $p$ in the other summand. These operations commute with $2$-sums along other
distinguished elements. The restriction of an ordering to a minor has no greater width.

Fix an edge $uv$, and write $E_u,E_v$ for the parts of the partition it displays.
Partition $E_u$ into groups indexed by the elements of $N_u$ other
than the distinguished element towards $v$. A remaining element gives a singleton group,
and any other distinguished element gives the set of remaining elements in the component beyond
that distinguished element. These groups are nonempty. Define the groups on $E_v$ in
the same way. There are $|E(N_u)|-1\geq3$ groups on one side and
$|E(N_v)|-1\geq3$ on the other. Choosing one representative in every
group and applying the preceding reduction gives a minor isomorphic
to $N_u\oplus_2N_v$, with peripheral distinguished elements renamed by the chosen
representatives. Proposition~\ref{prop:rooted-two-sum} forces the two
parts to be consecutive in every such restricted ordering.

This forces $E_u$ and $E_v$ to be consecutive as well. Otherwise,
after exchanging $u,v$ if necessary, there are elements $x,y,z$ in this
order with $x,z\in E_u$ and $y\in E_v$. If $x,z$ belong to different
groups, they and $y$ extend to a choice of representatives contradicting
the preceding paragraph. If they belong to the same group, take $w$
from a different group on $E_u$. When $w$ precedes $y$, use the triple
$w,y,z$, and otherwise use $x,y,w$. Again the two elements of $E_u$
belong to different groups, giving the same contradiction.
\end{proof}

\begin{theorem}[Tree decompositions using $3$-connected matroids]\label{thm:tree-two-sums}
The matroid $M(T)$ has path-width at most two exactly when $T$ is a
path and every matroid labelling a vertex has an ordering of width at most two with the distinguished elements for its incident edges at the ends. At an internal vertex, the two distinguished elements must occupy opposite ends.
For a one-vertex tree, the condition is simply $\pw(N_v)\leq2$.
Moreover, every ordering of $M(T)$ of width at most two lists the sets
$A_v$ consecutively in path order or its reverse.
\end{theorem}

\begin{proof}
Suppose that $M(T)$ has an ordering of width at most two. The one-vertex
case is immediate. At every leaf $v$, the set $A_v$ contains at least
three elements and is one part of the partition displayed by the incident edge. By
Lemma~\ref{lem:tree-shore-blocks}, this part is an initial segment or a final segment,
so it contains the first or the last element of the ordering. The sets $A_v$ for distinct leaves are disjoint. Thus $T$ has at most two leaves and is a
path, say $v_1,\ldots,v_t$.

Every $A_{v_i}$ is now nonempty, since at most two distinguished elements were removed from the matroid labelling $v_i$, which has at least four elements. Reverse the ordering if needed
so that $A_{v_1}$ is its first block. For each $i<t$, the set
$\bigcup_{j\leq i}A_{v_j}$ contains the first element and therefore
is an initial segment. Successive differences show that the sets $A_{v_i}$ are
consecutive in path order. Fix $i$ and reduce the rooted matroid obtained from the part of the path to its left, when present, to one representative, and do the same for the part to its right. The resulting minor is $N_{v_i}$ with its distinguished elements renamed.
Its inherited ordering has width at most two and consists of the left distinguished element, then $A_{v_i}$, then the right distinguished element, with absent distinguished elements omitted.
This proves the necessary endpoint conditions and the assertion about
all orderings.

Conversely, orient the stipulated orderings of the vertex labels along the path,
remove their distinguished elements, and concatenate the resulting orderings of the
sets $A_{v_i}$. Consider an initial segment ending within $A_{v_i}$, and let
$X_i$ be its intersection with that set. Let $L_i$ consist of the left
distinguished element of $N_{v_i}$ when $i>1$, and let $L_1=\varnothing$. The $2$-sum
rank formula gives
\begin{equation}\label{eq:tree-path-cut}
 \lambda_{M(T)}\!\left(\bigcup_{j<i}A_{v_j}\cup X_i\right)
 =\lambda_{N_{v_i}}(L_i\cup X_i).
\end{equation}
Indeed, the tuple $(d,a,b)$ in \eqref{eq:two-sum-cut} is $(0,1,0)$
when all elements of a rooted summand have been listed and $(0,0,1)$ when none have been listed. Substitution
places the distinguished element on the selected or unselected side, respectively.
Apply these identities successively to the rooted matroids on the left and right.
The right-hand side of \eqref{eq:tree-path-cut} is at most
two by the chosen ordering of $N_{v_i}$. Thus the concatenated ordering has
width at most two.
\end{proof}

The theorem also gives a local obstruction in more general
decompositions. If every vertex of a connected subtree in a decomposition using $2$-sums is labelled by a $3$-connected matroid on at least four elements, reduce the rooted matroids obtained from the components outside that subtree to representatives. This gives the sum over the subtree as a minor. Hence every such subtree must be a path when the original
matroid has path-width at most two.

The simultaneous requirement on two distinguished elements cannot be replaced by two
separate conditions allowing each distinguished element to occur at an end. Let $N=M(K_5\setminus04)$, with vertex set
$\{0,1,2,3,4\}$. Every element of $N$ can occupy an end of an ordering of width at most two. The edge symmetries have two orbits, represented by $01$ and
$12$, and the following orderings of width at most two start with these elements:
\[
\begin{split}
&01,02,03,12,13,23,14,24,34,\\
&12,01,02,03,13,23,14,24,34.
\end{split}
\]
Nevertheless, $\pw_{01,02}(N)>2$. To prove this, consider a three-element set of connectivity at most two. A dependent such set is a triangle. If it is independent, its rank is three, so its complement has rank at most three and the set contains a bond. Since $N$ is $3$-connected, a bond has at least three elements, and hence the set itself is a triad. The stars at $0$ and $4$ have three edges, the other vertex stars have four, and every cut with a two-vertex side has at least five edges. Thus the only triads are the stars at $0$ and $4$. If an ordering of width at most two
started with $01$ and ended with $02$, its first three elements would
therefore be $\{01,03,13\}$ and its last three would be
$\{02,03,23\}$. These triples share $03$, a contradiction.

Now form a $2$-sum of $N$ with a copy of $M(K_4)$ at each of $01$ and $02$. Both adjacent
$2$-sums have path-width two by
Proposition~\ref{prop:rooted-two-sum}, but their three-summand sum has
path-width greater than two by Theorem~\ref{thm:tree-two-sums}. For the upper bound, an ordering of $N$ with $01$ first and $02$ last and width three is
\[
01,03,13,12,23,14,24,34,02.
\]
Combining it with endpoint orderings of the two copies of $M(K_4)$ and using \eqref{eq:tree-path-cut} gives an ordering of the three-summand sum of width at most three. Its path-width is therefore exactly three. This
seventeen-element graphic example shows that requiring path-width at most two for the $2$-sum at each edge does not suffice even when the tree is a path.

The exclusion of circuit and cocircuit labels is essential. Start
with the cycle matroid of a triangle on $\{0,1,2\}$, and form its $2$-sums with copies of $M(K_4)$ at
$01,02,12$, using new vertex pairs $\{3,4\},\{5,6\},\{7,8\}$,
respectively. The resulting fifteen-element graphic matroid has a
three-leaf star as its decomposition tree, but it has path-width two. Supplementary Material S1 gives and
checks an ordering. Circuit and cocircuit summands therefore require further conditions for the path-width of their $2$-sums. Adler, Kant\'e and Kwon give algorithms for matroid path-width in
the circuit-and-cocircuit setting~\cite[Corollary~7.4]{AKK17}.
Their method uses a decomposition of the associated graphs into
complete graphs and stars and accounts for sets of elements obtained from tree components that are not consecutive in an ordering. Theorem~\ref{thm:tree-two-sums}
concerns the case in which every vertex is labelled by a $3$-connected matroid on at least four elements.

\subsection{Replacement of a three-connected summand}\label{sec:mixed-interfaces}

The preceding example shows that the elements of one summand can occur before and after elements of another summand in an ordering of width at most two. We now
consider a $2$-sum in which only one summand is required to be
$3$-connected. We first show when its elements can be made consecutive,
then use that ordering to compare its sums with all connected rooted matroids.

Throughout this subsection, $(N,p)$ is a rooted $3$-connected matroid with at least four elements. We form its $2$-sum with a connected rooted matroid $(B,p)$, taking the remaining ground sets disjoint.

\begin{lemma}[A consecutive ordering of a three-connected summand]\label{lem:rigid-block-normal}
Let $N$ be $3$-connected with at least four elements, and let $B$ be
a connected rooted matroid that is neither a circuit nor a cocircuit. If $N\oplus_2B$ has path-width at most
two, it has an ordering of width at most two in which $E(N)\setminus\{p\}$ is
consecutive. The elements of $B\setminus\{p\}$ then occupy at most
two blocks.
\end{lemma}

\begin{proof}
Put $A=E(N)\setminus\{p\}$. In an ordering of width at most two, move the first
element of $A$ immediately before the second. On every changed cut,
the triple for $N$ $(d,a,b)$ from \eqref{eq:two-sum-cut} changes from
$(1,0,1)$ to $(0,0,1)$, decreasing the cut value by one. Similarly,
move the last element of $A$ immediately after the penultimate one.
The triple changes from $(1,1,0)$ to $(0,1,0)$, again decreasing
the value by one. If $|A|=3$, the assertion follows.

Every remaining $B$-element between the first and last elements of $A$ now occurs at a cut
with at least two $A$-elements on each side. At a nontrivial split of
$B\setminus\{p\}$, connectedness gives $d_B\geq\max\{a_B,b_B\}$,
whereas $3$-connectivity gives
$d_N\geq1+\max\{a_N,b_N\}$. Substitution in
\eqref{eq:two-sum-cut} shows that width at most two permits only
\[
 (d_N,a_N,b_N)=(1,0,0),\quad(d_B,a_B,b_B)=(0,0,0),
\]
or
\[
 (d_N,a_N,b_N)=(2,1,1),\quad(d_B,a_B,b_B)=(1,1,1).
\]
Fix a $B$-element between the first and last elements of $A$, and consider the cuts
immediately before and after it. The split of $A$ is unchanged,
so the same one of the two displayed cases applies at both cuts
whenever the split of $B\setminus\{p\}$ is nontrivial. In the
first case, both cuts have connectivity zero in $B\setminus p$,
including an empty or complete cut when the element is first or
last among the $B$-elements. Since $B\setminus p$ is loopless,
the first zero cut makes adding this element increase the rank of the initial segment by one. The second zero cut then forces its removal to
decrease the rank of the final segment by one. Thus the element is a coloop of
$B\setminus p$, and hence in series with $p$ in $B$.
The dual argument shows that the second case makes the element
a loop of $B/p$, and hence in parallel with $p$.

Let $I$ be the subsequence of $B$-elements between the first and last elements of $A$, with each element of $I$ in series or in parallel with $p$ by the preceding argument. If $I$ is empty, the elements of $A$ are already consecutive. Otherwise either all elements of $I$ are in series with $p$ or all are in parallel with $p$. Indeed, a two-element circuit and a two-element cocircuit through $p$ cannot meet only in $p$. If the two pairs coincided, they would form a two-element component, contradicting the hypotheses that $B$ is connected and is neither a circuit nor a cocircuit. Let $P,Q$ be the $B$-subsequences before the first and after the last element
of $A$. At least one is nonempty, since otherwise $B$ would itself be
a circuit or a cocircuit.

If $P,Q$ are both nonempty, move $I$, in its original order, to the
end of $P$. The triple for $B$ at cuts between the first and last elements of $A$ remains $(0,0,0)$
in the series case or $(1,1,1)$ in the parallel case. Indeed, the
moved elements are coloops in the direct-sum separation of
$B\setminus p$, or loops in that of $B/p$, respectively. The new
cuts inside $I$ precede $A$ and have value one. If $P$ is nonempty
and $Q$ is empty, use the same move. Replacing the triple of the partial
$B$-ordering by the full triple changes each affected cut by
$-a_N$ in the series case and by $b_N-1$ in the parallel case.
If $P$ is empty, move $I$ to the beginning of $Q$, giving changes
$-b_N$ and $a_N-1$, respectively. These are nonpositive, and the
new cuts inside $I$ again have value one. The moves preserve the
relative order of all $B$-elements and make $A$ consecutive.
\end{proof}

To calculate the least width obtainable when $E(N)\setminus\{p\}$ is inserted as one consecutive block, define
\[
 h_N=\pw_p(N),\qquad d_N=\pw(N\setminus p),\qquad
 c_N=\pw(N/p).
\]
Here $d_N$ and $c_N$ denote path-widths of minors, not coordinates of a cut triple. For an ordering $\rho$ of $B$, let $w_B(\rho)$
denote its width. Let $(d_\rho,a_\rho,b_\rho)$ be the triple defined
in the proof of Proposition~\ref{prop:rooted-two-sum}, using the
elements before and after $p$ in $\rho$ and taking ranks in $B$.
The minimum possible largest cut value inside the inserted block is
given by the function
\[
 K_N(d,a,b)=
 \begin{cases}
 d+1+d_N,&(a,b)=(0,0),\\
 d+h_N,&(a,b)\in\{(0,1),(1,0)\},\\
 d+c_N,&(a,b)=(1,1).
 \end{cases}
\]
Write $N^{\mathrm s}_p$ and $N^{\mathrm p}_p$ for one series extension and one parallel extension at $p$, respectively. The following five indicators give the outcomes of the five conditions used for replacement, with one denoting a satisfied bound and zero a failed bound: \[
\begin{split}
 \omega_N&=[h_N\leq2],\qquad
 \delta_N=[d_N\leq1],\qquad \kappa_N=[c_N\leq1],\\
 \sigma_N&=[\pw(N^{\mathrm s}_p)\leq2],\qquad
 \pi_N=[\pw(N^{\mathrm p}_p)\leq2].
\end{split}
\]
We compare rooted matroids by the equivalence relation from the introduction, taking the ground sets outside the identified roots disjoint in each $2$-sum.

\begin{theorem}[Five conditions for a three-connected summand]\label{thm:rigid-interfaces}
Let $N$ be a rooted $3$-connected matroid with at least four elements.
If $B$ is connected and is neither a circuit nor a cocircuit, then
\[
 \pw(N\oplus_2B)\leq2
 \quad\Longleftrightarrow\quad
 \min_\rho\max\{w_B(\rho),
 K_N(d_\rho,a_\rho,b_\rho)\}\leq2.
\]
The five indicators determine the equivalence class of
every such $N$. Among rooted $3$-connected matroids with at least four
elements, there are at most eight classes.
\end{theorem}

\begin{proof}
For a consecutive block $E(N)\setminus\{p\}$, formula
\eqref{eq:two-sum-cut} gives its least possible maximum cut value as
$K_N(d,a,b)$. In the $(0,0)$ case this is one plus the deletion
connectivity, with the additive term $d$. In the $(1,1)$ case use
\[
 \lambda_{N/p}(X)=d'_N+1-a'_N-b'_N,
\]
where $(d'_N,a'_N,b'_N)$ is the triple of the split of
$E(N)\setminus\{p\}$. The two mixed cases give the endpoint width $h_N$,
with reversal exchanging the two ends. All cuts outside the block
are the corresponding cuts of $B$. The equivalence follows from
Lemma~\ref{lem:rigid-block-normal} and insertion of an optimal block.

For a nontrivial split of a connected $B$, one has
$d\geq\max\{a,b\}$. Deleting or contracting an element of $N$ gives a connected matroid with at least three elements, so $d_N,c_N\geq1$. Since $N$ is $3$-connected with at least four elements, every ordering has a cut of connectivity at least two, and hence $h_N\geq2$.
Consequently the displayed bound holds only when $p$ is at an end
and $\omega_N=1$, or when $(d,a,b)=(0,0,0)$ and $\delta_N=1$,
or when $(d,a,b)=(1,1,1)$ and $\kappa_N=1$. Thus $\omega_N,\delta_N,\kappa_N$ determine the outcome for every
connected $B$ that is neither a circuit nor a cocircuit. If $B$ is
a circuit with at least three elements, the $2$-sum replaces $p$
by a series class, and $\sigma_N$ determines the outcome. If $B$
is a cocircuit with at least three elements, $\pi_N$ does so.
The size of the series or parallel class beyond two does not change
the answer, by Theorem~\ref{thm:exact-doubling} and duality.
Finally, taking $B=U_{1,2}$ only renames $p$.

If $\pw(N)>2$, then $\pw(N\oplus_2B)>2$ for every connected
rooted $B$, because $N$ is a minor of the sum up to renaming $p$. Suppose
$\pw(N)\leq2$. Either $\delta_N=1$ or $\kappa_N=1$ implies
$\omega_N=1$, by placing $p$ at an end of an ordering of the deletion
or contraction. In turn, $\omega_N=1$ implies
$\sigma_N=\pi_N=1$, since an endpoint may be expanded without
increasing width beyond two.

At least one of $\sigma_N,\pi_N$ is one. Indeed, in an ordering of width at most two of $N$, the triple $(d,a,b)$ at $p$ satisfies
$d\leq1+\min\{a,b\}$ from the cuts immediately before and after
$p$. If $d\leq1$, adding an element in series with $p$ gives an internal cut value $d+1\leq2$. Otherwise $d=2$ and $a=b=1$, so adding an element in parallel with $p$ gives an internal cut value $d+2-a-b=2$.

There are therefore at most four possible combinations of indicators with $\omega_N=1$, indexed by $(\delta_N,\kappa_N)$, and at most three with $\omega_N=0$, indexed by the nonzero pair $(\sigma_N,\pi_N)$. Together with the case $\pw(N)>2$, this gives the bound eight. No assertion that every listed combination is attained is needed.
\end{proof}

For a $3$-connected summand on at least four elements, the theorem
reduces replacement in an arbitrary connected $2$-sum to five explicit conditions. General trajectory methods already retain finite information
for replacement~\cite{JKO17}. Here the consecutive ordering permits
the simpler description by deletion, contraction, an endpoint
condition, series extension and parallel extension. The theorem makes the elements of one summand consecutive, without asserting that
the elements from all the other components of a decomposition tree can be made consecutive in the same ordering.

The next example shows why deletion, contraction and endpoint widths
alone do not determine replacement for arbitrary connected rooted matroids. Let $B_1=M(K_4)$, rooted at $01$, and let
$B_2=M(K_4-23)$ with the same root. Both have
\[
 (\pw(B_i\setminus p),\pw(B_i/p),\pw_p(B_i))=(1,1,2).
\]
For $B_1\setminus p$, the order $02,03,23,12,13$ has width one.
For $B_1/p$, list the two parallel pairs on opposite sides of $23$.
For $B_2$, deletion gives a circuit and contraction gives two parallel
pairs. Reverse any ordering of $B_2$ that begins with $p$, and take the cut after the first three elements of the reversed ordering. Those three edges have rank three and their complement has rank two, so the cut has value two and the endpoint width is at least two. It is at most two because $B_2$ is a restriction of
$B_1$ and an endpoint ordering of $B_1$ restricts to one of $B_2$.

Form a connected rooted matroid $C$ by taking the $2$-sums of a central copy of $M(K_4)$ with two further copies at edges $01$ and $02$, using new vertex pairs $\{4,5\}$ and $\{6,7\}$, respectively, and retaining $p=03$ as the root. The resulting decomposition tree is a path with three vertices. Form $C\oplus_2B_2$ at $p$, using new vertices $8,9$. The resulting seventeen-element matroid has the ordering of width at most two \[
\begin{split}
67,27,26,07,06,39,09,23,12,13,38,08,15,14,45,05,04.
\end{split}
\] The successive cut values are $1,2,2,2,1,2,2,2,2,2,2,1,2,2,2,1,0$.
Replacing $B_2$ by $B_1$ adds edge $89$. The decomposition tree is
then a three-leaf star of four copies of $M(K_4)$, so the resulting matroid has path-width greater than two by Theorem~\ref{thm:tree-two-sums}. Adding the one edge
to the displayed ordering gives an upper bound of three.

The two rooted matroids therefore have the same deletion, contraction and
endpoint parameters, yet their sums with $C$ have different answers at
width two. In the displayed ordering, the elements of
$B_2=M(K_4-23)$ occur on both sides of elements of the central summand.
If those elements could instead form one consecutive block, its
minimum cost would depend only on the three equal parameters.
Replacing that block by one from $B_1$ would then give an ordering
of width at most two, a contradiction. In the next section we therefore retain the cut information needed
to determine how the orderings of the two summands can be interleaved.

\section{Finite tests and a recurrence for two-sums}
\label{sec:finite-interfaces}

Let $(N,p)$ be a connected rooted matroid, where $p$ is neither a
loop nor a coloop. Two such matroids are \emph{equivalent}
if
\[
 \pw(N\oplus_2 C)\leq2
 \quad\Longleftrightarrow\quad
 \pw(N'\oplus_2 C)\leq2
\]
for every connected rooted matroid $(C,p)$. The ground sets outside
$p$ are taken disjoint when forming a $2$-sum. Thus equivalence means that replacing one rooted matroid by the other in any such $2$-sum does not change whether the path-width is at most two. Our aim is to
decide this universal condition by finitely many specified tests.

We give an explicit finite test for this equivalence and an exact
recurrence for arbitrary tree decompositions using $2$-sums. For each ordering, we retain a sequence of rank and closure data at
its cuts. The compression follows the trajectory method of Jeong,
Kim and Oum~\cite{JKO17}, here expressed by rank and closure data for a $2$-sum at one root.
Kant\'e, Kim, Kwon and Oum compare sets of compact trajectories in their replacement argument~\cite[Lemma~3.4]{KKKO23}. Here we compute the
resulting tests explicitly at width two. There are thirteen fixed
graphic or cographic rooted matroids, each on at most fourteen elements,
and the matroid being tested need not be representable.

For $A=E(N)\setminus\{p\}$ and $X\subseteq A$, put $Y=A\setminus X$
and define
\[
 d=r_N(X)+r_N(Y)-r(N),\qquad
 a=\mathbf1_{p\in\operatorname{cl}_N(X)},\qquad
 b=\mathbf1_{p\in\operatorname{cl}_N(Y)}.
\]
For an ordering of $A$, form the sequence of triples $(d,a,b)$ at every initial segment, including the empty and complete initial segments, to test how the ordering can be interleaved with an ordering of another summand. The sequence starts at
$(0,0,1)$ and ends at $(0,1,0)$. The coordinate $a$ is nondecreasing
and $b$ is nonincreasing. For a nonempty proper initial segment, connectedness
gives
\begin{equation}\label{eq:interface-connected}
 d\geq\max\{a,b\},
\end{equation}
since $\lambda_N(X)=d+1-b$ and
$\lambda_N(X\cup\{p\})=d+1-a$ are positive.

For triples $u=(d,a,b)$ and $v=(e,c,f)$, the $2$-sum rank formula gives the connectivity of the corresponding cut as \begin{equation}\label{eq:interface-cell}
 \Phi(u,v)=d+e+1-ac-bf.
\end{equation} To interleave two orderings, choose at each step the next element of
one ordering, preserving the order within each summand. If the cut
sequences are $u_0,\ldots,u_s$ and $v_0,\ldots,v_t$, the current
position is a pair $(i,j)$ in the rectangular grid
$\{0,\ldots,s\}\times\{0,\ldots,t\}$. Appending an element advances
one coordinate by one, and the cut value is $\Phi(u_i,v_j)$.
The two sequences admit an interleaving of width at most two precisely when the grid has a path from $(0,0)$ to $(s,t)$, advancing one coordinate by one at each step and using only positions with $\Phi\leq2$.

If two sequences admit such an interleaving, every triple has $d\leq2$. Indeed, the contribution
$e+1-ac-bf$ in \eqref{eq:interface-cell} is nonnegative by
\eqref{eq:interface-connected}, with the endpoint cases checked
directly. A triple with $(a,b)=(0,0)$ must further satisfy $d\leq1$.
We therefore restrict to sequences satisfying these bounds at every position. Sequences arising from orderings also satisfy the endpoint and monotonicity conditions established above.

\begin{lemma}[Compression at width two]\label{lem:interface-compression}
Every sequence arising from an ordering and satisfying the stated bounds can be reduced to at most nine triples without changing which sequences admit an interleaving of width at most two with it. All the reduced sequences lie in an explicit set of $114$ possibilities, independent of the matroid.
\end{lemma}

\begin{proof}
Within a consecutive interval having constant $(a,b)$, delete
consecutive repetitions of $d$. Also delete the entries strictly
between positions $i<j$ whenever all those entries lie between
$d_i$ and $d_j$.

These deletions preserve the existence of an interleaving of width at most two with every fixed external sequence. If $d_i\leq d_j$, start with a feasible interleaving of the original interval, perform all external advances while the local value of $d$ is $d_i$, then advance to $d_j$ at the original final external position.
All substituted cell values are no larger. Conversely, in an interleaving
of the reduced interval insert the deleted triples immediately before
the advance to $d_j$, holding the external position fixed. Their
values are at most $d_j$, so these cells are feasible. For
$d_i\geq d_j$, perform the external advances at $d_j$ and reinsert
the deleted triples at the external position of $d_i$. Repetitions
are handled in the same way. The argument also shows how to expand a successful reduced interleaving into an ordering of the original elements.

The monotonicity of $a$ and $b$ gives at most three consecutive
intervals on which $(a,b)$ is constant. Write $01$, for example, for $(a,b)=(0,1)$. Their order
is $01,10$, or $01,00,10$, or $01,11,10$. In the initial interval the first
$d$ is zero and subsequent values are one or two. Its irreducible
sequences are
\[
 0,\quad 01,\quad 02,\quad 021.
\]
The final interval has the four reverse sequences. The middle interval uses
$\{0,1\}$ in case $00$ or $\{1,2\}$ in case $11$. On either
two-element set $\{x,y\}$ its irreducible sequences are
\[
 x,\quad y,\quad xy,\quad yx,\quad xyx,\quad yxy.
\]
Hence every reduced sequence has length at most nine, and
$4\cdot4\cdot(1+6+6)=208$ sequences form a containing set.

At a single-element advance where $a$ changes from zero to one, the
rank of the initial segment increases by one. Where $b$ changes from one to zero, the
rank of the final segment decreases by one. Consequently, across the boundaries
of intervals the possible changes in $d$ are
\[
\begin{array}{c|ccccc}
\text{transition}&01\to10&01\to00&00\to10&01\to11&11\to10\\
\hline
\Delta d&0&-1,0&0,1&0,1&-1,0.
\end{array}
\]
Compression retains the triples at the ends of each interval, so the reduced sequences must satisfy these same transition conditions. Exactly $114$ of the $208$ sequences in the containing set do so. For completeness,
the multiplicities of terminal values $0,1,2$ in the initial interval
are $1,2,1$, as are the multiplicities of initial values in the final
interval. The direct transition contributes $1^2+2^2+1^2=6$ sequences.
Each of the six middle sequences in either middle interval has three
allowable initial intervals and three allowable final intervals, giving
$6+6\cdot9+6\cdot9=114$.
\end{proof}

Let $\mathcal W$ denote this explicit set of $114$ sequences.
Reversing an ordering reverses its sequence of triples and exchanges
the coordinates $a$ and $b$. To compare families of reduced sequences, define, for $\mathcal F\subseteq\mathcal W$, the set of sequences that can be paired with at least one member of the family by
\[
 \mathcal L(\mathcal F)
 =\{v\in\mathcal W: v\text{ admits such an interleaving with some }u\in\mathcal F\}.
\]
For a rooted matroid $(N,p)$, let $\mathcal F$ consist of the reduced sequences obtained from all orderings of $E(N)\setminus\{p\}$ whose cut sequences satisfy the stated bounds. Write the resulting set as $\mathcal L(N,p)$. Each reduction preserves the existence of such interleavings, so its order does not change this set. Reversing an
ordering gives another ordering of the same ground set, and hence
the family is closed under reversal.

We first specify the matroids that will test equivalence at width two. In Table~\ref{tab:interface-tests}, the pair $(G,e)$ denotes
the cycle matroid of the displayed graph rooted at edge $e$. A star
on a rooted matroid denotes its dual, with the same root. All graphs
are simple. Define $C$ by starting with $K_4$ on $\{0,1,2,3\}$ and
attaching a copy of $K_4$ at $01$ with new vertices $4,5$ and another
at $02$ with new vertices $6,7$, deleting the common edge used for each $2$-sum. Its edges
are
\[
 12,13,23,04,05,14,15,45,06,07,26,27,67,03.
\]
Define $D$ to have the edge $01$, all edges from $\{0,2\}$ to
$\{3,4,5\}$, and all edges from $\{2,1\}$ to $\{6,7,8\}$.
Thus $D$ consists of two copies of $K_{2,3}$ sharing vertex $2$, together with edge $01$, which is the root in the corresponding cycle matroid.

\begin{table}[t]
\centering\small
\begin{tabular}{ccl}
\hline
\text{Matroids}&\text{Size}&\text{Rooted matroids}\\
\hline
$T_1,T_2$&12&$(P_6^3,14)$ and its dual\\
$T_3,T_4$&14&$(C,03)^*$ and $(C,03)$\\
$T_5,T_6$&13&$(D,01)$ and its dual\\
$T_7,T_8$&7&$(K_{2,3}+01,01)^*$ and $(K_{2,3}+01,01)$\\
$T_9,T_{10}$&5&$(K_4\setminus23,01)^*$ and $(K_4\setminus23,01)$\\
$T_{11}$&5&$(K_4\setminus23,02)$\\
$T_{12},T_{13}$&3&$M(C_3)$ and $M(C_3)^*$, rooted at any element\\
\hline
\end{tabular}
\caption{The thirteen fixed matroids. The graph $P_6^3$ has vertex set
$\{0,\ldots,5\}$ and edges $ij$ with $1\leq j-i\leq3$.
The two parts of $K_{2,3}$ are $\{0,1\}$ and $\{2,3,4\}$.}
\label{tab:interface-tests}
\end{table}

\begin{samepage}
\begin{theorem}[Thirteen fixed matroids for equivalence]
\label{thm:finite-interface}
For connected rooted matroids $(N,p)$ and $(N',p)$, the following are
equivalent:
\begin{enumerate}
\item They are equivalent.
\item For every $1\leq j\leq13$,
$\pw(N\oplus_2T_j)\leq2$ if and only if
$\pw(N'\oplus_2T_j)\leq2$.
\item $\mathcal L(N,p)=\mathcal L(N',p)$.
\end{enumerate}
There are exactly twenty equivalence classes, each represented by a
graphic or cographic matroid.
\end{theorem}
\end{samepage}

\begin{proof}
By Lemma~\ref{lem:interface-compression}, the sets $\mathcal L(N,p)$ and $\mathcal L(N',p)$ determine whether each of the corresponding $2$-sums has path-width at most two. Thus the third assertion implies the first, which implies the second.

We give the finite reduction proving the remaining implication.
It involves only sequences of length at most nine and grids of at
most $81$ positions. For each $u\in\mathcal W$, calculate the set of
$v\in\mathcal W$ admitting an interleaving of width at most two with $u$, using
\eqref{eq:interface-cell}. There are $22$ distinct such sets.
Next form every union permitted by taking a family of sequences
together with their reversals. This gives precisely the $22$ sets
listed in Table~\ref{tab:interface-codes}, including the empty set.
In each row, the listed sequences and their reversals
generate that union of sets. The vector of outcomes for the thirteen fixed matroids has entry one in position $j$ exactly when one of those sequences admits an interleaving of width at most two with a sequence from $T_j$.
For a row realised by $(N,p)$, this says precisely that
$\pw(N\oplus_2T_j)\leq2$.
These finite calculations are supplied with independently implemented
checks in Supplementary Material~S1. They enumerate the displayed
finite sequence set and compute the subset ranks of the thirteen fixed matroids, so the coverage assertion does not depend on a bounded
search through input matroids.

Two rows of the table cannot arise from a connected rooted matroid.
Tutte's connected deletion or contraction theorem states that at least
one of $N\setminus p$ and $N/p$ is connected \cite{Oxl84}. If
$N\setminus p$ is connected, a cut sequence has no interior triple
$(0,0,0)$. Taking reversal-closed unions of rows of sequences avoiding that
triple gives the sets with indices
\[
 0,1,2,3,4,5,6,7,9,11,12,13,14,16,18,20,21.
\]
If $N/p$ is connected, it has no interior triple $(1,1,1)$, because
$\lambda_{N/p}(X)=d+1-a-b$. The corresponding possible indices are
\[
 0,1,2,3,4,5,6,7,8,11,12,13,14,15,18,19,21.
\]
The union of these two lists excludes exactly indices $10$ and $17$.
The remaining twenty outcome vectors in Table~\ref{tab:interface-codes}
are distinct. Consequently the second assertion determines the whole
set $\mathcal L(N,p)$ and implies the third.

Each of the twenty remaining sets has an explicit graphic or
cographic representative. Nineteen representatives are listed below,
and Proposition~\ref{prop:class-three} supplies the remaining one.
Their complete families of cut sequences and pairwise distinguishing tests are
verified in S1. Since the thirteen tests distinguish their sets,
these representatives belong to twenty distinct equivalence classes.
\end{proof}

For the finite table, encode the ten possible triples by the following symbols:
\[
\begin{array}{c|cccccccccc}
\text{symbol}&i&f&l&L&r&R&s&S&p&P\\
\hline
(d,a,b)&001&010&101&201&110&210&000&100&111&211.
\end{array}
\]
Here each three-digit entry denotes the corresponding triple, rather
than an integer. Reversal exchanges $i,f$, $l,r$, and $L,R$, fixes
the other symbols, and reverses their order.

\begin{table}[t]
\centering\small
\begin{tabular}{ccl}
\hline
\text{Index}&\text{Listed sequences}&\text{Outcomes for $T_1,\ldots,T_{13}$}\\
\hline
$0$&$\varnothing$&$\mathtt{0000000000000}$\\
$1$&$\mathtt{iLSRf}$&$\mathtt{0000000000010}$\\
$2$&$\mathtt{iLPRf}$&$\mathtt{0000000000001}$\\
$3$&$\mathtt{iLSRf},\ \mathtt{iLPRf}$&$\mathtt{0000000000011}$\\
$4$&$\mathtt{iLSrRf}$&$\mathtt{0000000000111}$\\
$5$&$\mathtt{iLlSrRf}$&$\mathtt{0000000010111}$\\
$6$&$\mathtt{iLlPrRf}$&$\mathtt{0000000001111}$\\
$7$&$\mathtt{iLlrRf}$&$\mathtt{0000000011111}$\\
$8$&$\mathtt{iLSsSRf}$&$\mathtt{0000001011111}$\\
$9$&$\mathtt{iLPpPRf}$&$\mathtt{0000000111111}$\\
$10^{\dagger}$&$\mathtt{iLSsSRf},\ \mathtt{iLPpPRf}$&$\mathtt{0000001111111}$\\
$11$&$\mathtt{isSRf}$&$\mathtt{0000001111111}$\\
$12$&$\mathtt{isSsf}$&$\mathtt{0000101111111}$\\
$13$&$\mathtt{ipPpf}$&$\mathtt{0000011111111}$\\
$14$&$\mathtt{isSsf},\ \mathtt{ipPpf}$&$\mathtt{0000111111111}$\\
$15$&$\mathtt{ilsrf}$&$\mathtt{0010111111111}$\\
$16$&$\mathtt{ilprf}$&$\mathtt{0001111111111}$\\
$17^{\dagger}$&$\mathtt{ilsrf},\ \mathtt{ilprf}$&$\mathtt{0011111111111}$\\
$18$&$\mathtt{isrf}$&$\mathtt{0011111111111}$\\
$19$&$\mathtt{isf}$&$\mathtt{1011111111111}$\\
$20$&$\mathtt{ipf}$&$\mathtt{0111111111111}$\\
$21$&$\mathtt{if}$&$\mathtt{1111111111111}$\\
\hline
\end{tabular}
\caption{The complete finite reduction. Each row gives the union of the sets defined by its listed sequences and their reversals. The daggered rows cannot occur for a connected rooted matroid. The twenty undaggered outcome vectors are distinct. The listed sequences are formal and are not asserted to be cut sequences of one matroid.}
\label{tab:interface-codes}
\end{table}

Nineteen of the rows are witnessed by the thirteen fixed rooted matroids and six additional rooted matroids. The matroids $M(K_5)$, $M(K_4)$ and $U_{1,2}$, rooted at any element, realise rows $0,14,21$, respectively. Taking $13$ as the root of the cycle matroid of $C$ gives a representative for row $4$. For row $7$, construct $C'$ as $C$, using the disjoint edges $01,23$ instead of $01,02$ as the common edges for the $2$-sums, and take its cycle matroid with root $03$. Finally, row $11$ is realised by $(E,01)$, where $E$ is $K_4$ on $\{0,2,3,4\}$ together with the two-edge path $0,1,2$.

The remaining row has a small graphic representative. Start with $K_4$
on $\{0,1,2,3\}$, retain all six edges, and add a new two-edge path
alongside each of $01,02,03,23$, with new vertices $4,5,6,7$, respectively.
Write $G_3$ for the resulting graph, so that
\[
 E(G_3)=\{01,02,03,12,13,23,04,14,05,25,06,36,27,37\}.
\]

\begin{proposition}[A graphic representative for row three]
\label{prop:class-three}
The connected rooted matroid $(N,p)=(M(G_3),12)$ has fourteen elements
and rank seven, and realizes row $3$ of Table~\ref{tab:interface-codes}.
In particular,
\[
 \pw(N_p^{\mathrm s})=\pw(N_p^{\mathrm p})=2,
 \qquad \pw(N\oplus_2T_{11})=3.
\]
\end{proposition}

\begin{proof}
Adding paths between existing vertices preserves $2$-connectedness,
so $N$ is connected and has the stated size and rank. We give the
finite calculation identifying its equivalence class. Compute every subset rank from
$r(X)=|V(G_3)|-c(V(G_3),X)$, where $c(V(G_3),X)$ is the number
of connected components of the spanning subgraph with edge set $X$,
including isolated vertices. The complete family of reduced sequences satisfying the stated bounds is
\[
 \{\mathtt{iLSRf},\mathtt{iLPRf},\mathtt{iLRf}\}.
\]
The first two sequences are listed in row $3$. The third adds no external sequence admitting an interleaving of width at most two. Indeed, the triple $L=(2,0,1)$ can be
paired only with the external initial triple $i$, and $R=(2,1,0)$ only
with the external final triple $f$. Since $L$ and $R$ are consecutive,
no grid path satisfying $\Phi\leq2$ can move between those pairings. Thus the
computed family gives precisely row $3$.

For a separate check of the three displayed widths, the series graph
replaces $12$ by $18,82$, the parallel graph adds a second copy of $12$,
and the sum with $T_{11}$ replaces $12$ by $18,19,89,82$. In each graph, compute all subset ranks by component counts.
Let $W(X)$ be the minimum, over orderings of $X$, of the largest
connectivity value of an initial segment, with connectivity evaluated in the
full matroid. Then use
\[
 W(\varnothing)=0,\qquad
 W(X)=\max\bigl\{\lambda(X),\min_{e\in X}W(X\setminus\{e\})\bigr\},
\]
with every cut evaluated in the full matroid. Choosing the last element
of an ordering proves this recurrence by induction on $|X|$, so its
final value is the exact path-width. The matroid $N$, its series extension and its parallel extension all have path-width two, while $N\oplus_2T_{11}$ has path-width three. The standalone program
\texttt{check\_class3\_witness.py} in S1 performs both calculations and
supplies optimal orderings and all initial segment cut values. The verifier for the complete finite table independently computes the same complete family of sequences by
binary elimination and checks its set against all $114$ sequences.
\end{proof}

The finite classification also determines the minimum number of fixed
tests. Each such test asks whether the $2$-sum with one fixed connected
rooted matroid has path-width at most two.

\begin{corollary}[Optimal number of fixed tests]
\label{cor:minimum-tests}
Any family of fixed connected rooted matroids whose $2$-sums distinguish all equivalence classes contains at least thirteen members.
Thus the thirteen tests in Table~\ref{tab:interface-tests} form a family of
minimum cardinality.
\end{corollary}

\begin{proof}
Let $I=\{0,1,\ldots,21\}\setminus\{10,17\}$ be the indices of the
realizable rows of Table~\ref{tab:interface-codes}. By
Theorem~\ref{thm:finite-interface}, every connected rooted matroid used for a fixed test
belongs to one of these classes. Moreover, two tests in the same class
give the same answer on every connected rooted input, since $2$-sum
is commutative. It therefore suffices to prove that any set of test
classes distinguishing all twenty input classes has at least thirteen
members.

For distinct $a,b\in I$, let $D(a,b)\subseteq I$ consist of the test
classes that give different answers on input classes $a$ and $b$.
These sets are calculated from the listed sequences in
Table~\ref{tab:interface-codes} and their reversals, using
\eqref{eq:interface-cell}. A test class gives answer one on an input
class exactly when some listed sequence of one admits an interleaving of width at most two with some listed sequence of the other. The calculation gives
\[
\begin{array}{lll}
D(19,21)=\{2\},&D(20,21)=\{1\},&D(11,12)=\{8\},\\
D(11,13)=\{9\},&D(4,5)=\{15\},&D(4,6)=\{16\},\\
D(3,4)=\{18\},&D(2,3)=\{19\},&D(1,3)=\{20\}.
\end{array}
\]
To distinguish each of these pairs, the family must therefore contain
a test from each of the nine classes
\[
\{1,2,8,9,15,16,18,19,20\}.
\]
The same calculation gives
\[
\begin{array}{lll}
D(16,18)=\{4,5\},&D(15,18)=\{4,6\},&D(15,16)=\{5,6\},\\
D(9,11)=\{11,12\},&D(8,11)=\{11,13\},&D(8,9)=\{12,13\}.
\end{array}
\]
The first line requires tests from at least two of the classes
$\{4,5,6\}$, since one class cannot belong to all three displayed
pairs. Similarly, the second line requires at least two of
$\{11,12,13\}$. These two sets are disjoint from one another and from
the nine compulsory classes, so every distinguishing family has at
least $9+2+2=13$ members. Theorem~\ref{thm:finite-interface} shows that
the thirteen tests in Table~\ref{tab:interface-tests} attain this bound.
\end{proof}

We now use the reduced sequences to combine an arbitrary tree decomposition using $2$-sums. The recurrence needs the positions within each sequence,
so it retains families of sequences, rather than only the corresponding sets $\mathcal L$.

This applies in particular to the canonical tree decomposition of a connected matroid with at least three elements, whose vertices are labelled by circuits, cocircuits and $3$-connected matroids on at least four elements. Each matroid labelling a vertex has at least three elements, no two vertices labelled by circuits are adjacent, and no two vertices labelled by cocircuits are adjacent. With these conditions the decomposition is unique up to the usual relabelling of the distinguished elements~\cite{Oxl11}. A connected matroid on one or two elements may be represented by a one-vertex tree, and the two-element case is $U_{1,2}$. The recurrence also permits redundant vertices labelled by $U_{1,2}$ in a decomposition using $2$-sums.

Choose a root vertex of the tree $T$. For a vertex $v$ other than the root vertex, let $H$ be its matroid label. Let $q$ be the distinguished element for the edge to its parent and let $p_1,\ldots,p_t$ be the distinct distinguished elements for the edges to its children. Put $R=E(H)\setminus\{q,p_1,\ldots,p_t\}$. For each $j$, let $B_j$ be the connected rooted matroid obtained by taking the $2$-sums of the labels in the corresponding child subtree, retaining $p_j$ as its root. To describe a cut after
combining these matroids, choose the listed subset $X\subseteq R$
and a position in each child sequence, whose triple is
$(d_j,a_j,b_j)$. Put
\[
 Y=R\setminus X,\qquad
 A=\{p_j:a_j=1\},\qquad B=\{p_j:b_j=1\}.
\]
Here $A$ and $B$ contain the distinguished elements in the closure of
the listed and unlisted child elements, respectively. The triple
for the combined rooted matroid is
\begin{equation}\label{eq:finite-factor-state}
\begin{split}
 d_H={}&r_H(X\cup A)+r_H(Y\cup B)-r(H)
       +\sum_{j=1}^t(d_j+1-a_j-b_j),\\
 a_H={}&\mathbf1_{q\in\operatorname{cl}_H(X\cup A)},\qquad
 b_H=\mathbf1_{q\in\operatorname{cl}_H(Y\cup B)}.
\end{split}
\end{equation}
At the root vertex there is no distinguished element for a parent edge, so the first line is the connectivity value and the other two coordinates are omitted.

\begin{theorem}[Recurrence for tree decompositions using two-sums]
\label{thm:general-two-sum-tree}
Let a connected matroid be given by a tree decomposition using $2$-sums of connected matroids, with redundant vertices labelled by $U_{1,2}$ also permitted. This includes its canonical tree decomposition. The following procedure accepts exactly when the matroid has path-width at most two and constructs an ordering of width at most two when it accepts.

Process the vertices of the rooted tree from the leaves towards the root vertex. At a vertex other than the root vertex, choose one sequence from the family associated with each child subtree. Start with no element of $R$ listed and each child sequence
at its first position. A move either adds one element of $R$ to
the listed subset or advances one child sequence by one position.
These choices form the product of the subset lattice of $R$ and
the chains of positions in the chosen sequences.

Assign the triple in \eqref{eq:finite-factor-state} to every position
of this product. Discard positions violating $d_H\leq2$, or
$d_H\leq1$ when $(a_H,b_H)=(0,0)$. For every path from the initial
position to the position where all elements and child sequences
are complete, retain the reduced sequence of its triples. Repeat
this for every choice of child sequences to obtain the family
for the subtree at this vertex. At the root vertex, which has no parent element, use the first line of \eqref{eq:finite-factor-state} as the cut value and accept exactly when a path from the initial to the final position has value at most two throughout.
Moreover, the rooted matroid represented by any subtree below a nonroot vertex may be replaced by another connected rooted matroid with the same vector of outcomes for the thirteen fixed matroids without changing acceptance.
\end{theorem}

\begin{proof}
The rank of a global initial segment is
\[
 r_H(X\cup A)+
 \sum_{j=1}^t\bigl(r_{B_j}(X_j)-a_j\bigr),
\]
where $X_j$ is its intersection with child $j$. The corresponding
formula holds for the final segment, and the rank of the matroid obtained from this subtree is
$r(H)+\sum_j(r(B_j)-1)$. Adding the ranks of the initial and final segments and subtracting the rank of this matroid gives the first line of \eqref{eq:finite-factor-state}. Adding $q$ to the initial segment or final segment changes only the corresponding central rank term, which gives $a_H$ and $b_H$.

Before compression, every path in the product is exactly an interleaving
of the chosen child orderings and an ordering of the remaining elements.
Conversely, every ordering restricts to such a path. After child sequences have been compressed, a single advance can
pass over several elements. The induction must therefore compare
the stored family with all sequences that could be supplied by
the rest of the sum. Such an external sequence begins at $(0,0,1)$,
ends at $(0,1,0)$, has nondecreasing $a$ and nonincreasing $b$, and
satisfies \eqref{eq:interface-connected} at each interior position.
Consecutive triples need not arise by adding a single element of
one matroid. The induction asserts that an external sequence admits an interleaving of width at most two with a stored sequence if and only if it admits such an interleaving with a sequence from an actual ordering of the subtree matroid. We may restrict to external sequences satisfying the stated bounds, since every triple in a sequence admitting such an interleaving satisfies those bounds. Each permitted triple occurs at a cut of one of the connected rooted matroids $U_{1,2}$, $U_{2,3}$, $U_{1,3}$ and $M(K_4)$, although the triples need not occur together in a single ordering. Their separate realizations justify
the inequalities following from connectedness when the rank formula combines them.

To prove this invariant, fix an external sequence and an order of moves for all
the other summands. Applying \eqref{eq:finite-factor-state} to those fixed choices gives a sequence with the stated endpoints and monotone closure coordinates, satisfying \eqref{eq:interface-connected} at every interior position.
Regard it as the external sequence for the remaining child. By the inductive invariant,
a feasible interleaving using a stored child sequence exists if and only if
one using an actual child ordering exists, with the same order of all
external advances. Perform these substitutions successively for the
children, then compress the resulting parent sequence. The two-direction
argument of Lemma~\ref{lem:interface-compression} preserves the existence of an interleaving of width at most two with the fixed external sequence. This proves the invariant
inductively.
Every stored triple is realized by an actual subset of the matroid obtained from the subtree. Indeed, each child triple can separately be realized at a
subset of that child, and the rank formula then realizes the parent
triple. These subset witnesses need not form one ordering, but the
pointwise property justifies discarding triples that violate the bounds following from connectedness. The external contribution in \eqref{eq:interface-cell} is
nonnegative for every such external triple, and is at least one when
the local closure indicators satisfy $(a,b)=(0,0)$. Thus no
interleaving of width at most two can use a discarded triple.

At the root vertex there is no external summand, so the induction
assertion gives exactly the existence of an ordering of width at
most two. To construct that ordering, retain the local paths before
reduction and the choices used for each child. Expand each compressed child sequence by the feasible grid path from Lemma~\ref{lem:interface-compression}, preserving the order of all elements outside that child. Repeating the expansion recovers
an ordering of the ground set of the whole matroid. Finally, the labels outside a subtree below a nonroot vertex also give a connected rooted matroid, with the distinguished element for the edge to that subtree retained as its root.
Theorem~\ref{thm:finite-interface} therefore proves the replacement
assertion.
\end{proof}

Every stored reduced family lies in the fixed containing set of $208$ sequences from Lemma~\ref{lem:interface-compression}. A compressed child advance may span several elements, so an intermediate sequence need not satisfy the single-element conditions at a change of $(a,b)$. Those conditions reduce the containing set of $208$ sequences to the $114$ possibilities for actual orderings. Thus the recurrence
stores a family drawn from a fixed finite set, independently of
the number of elements in the matroid obtained from a subtree. The product of the subset lattice and sequence-position chains still depends on the size of the matroid labelling the vertex and the degree of that vertex. The theorem asserts an exact characterization and a constructive recurrence, without a claim that this product has polynomial size when those quantities are unrestricted.

The finite-test theorem also bounds the diameter of the decomposition
tree of an excluded minor. We first explain how to remove the middle of a
decomposition path while retaining the $2$-sums with the rooted matroids obtained from its two outer components.
This uses the replacement strategy for matroid path-width in
\cite[Lemma~3.4 and Proposition~5.1]{KKKO23}. Bounds on nested
separations also control decomposition paths for representability
obstructions in \cite{BDG16}. Here the explicit table for width two gives
a small bound without a representability assumption.

\begin{lemma}[Shortening a decomposition path]
\label{lem:interface-splice}
Let $T$ be the tree of a decomposition using $2$-sums that represents a matroid $M$, with each vertex labelled by a connected matroid on at least three elements and with distinct distinguished elements for its incident edges.
Orient a path with edges $e_1,\ldots,e_m$, and let $A_i,B_i$ be
the rooted matroids obtained from the components of $T-e_i$ on the left and right of $e_i$, respectively, retaining the distinguished element for $e_i$ as their common root. If $i<j$, then $A_i\oplus_2B_j$ is a proper
minor of $M$, after identifying the roots. Moreover, $A_i$ is a
proper minor of $A_j$ obtained without deleting or contracting its root, and $B_j$ is such a proper minor of $B_i$, up to relabelling the roots.
\end{lemma}

\begin{proof}
Let $H$ be the sum over the component between $e_i$ and $e_j$,
with its two outer distinguished elements $p_i,p_j$ retained. It is connected, and
$M=A_i\oplus_{2,p_i}H\oplus_{2,p_j}B_j$. If the vertices in that component are indexed by $S$, then
\[
 |E(H)|=2+\sum_{v\in S}(|E(N_v)|-2)\geq3.
\]
A circuit $Q$ of $H$ contains both distinguished elements. Delete $E(H)\setminus Q$
and contract $Q\setminus\{p_i,p_j\}$ to leave $U_{1,2}$ on the
distinguished elements. These operations involve only elements of $H$ other than $p_i,p_j$, all of which belong to $E(M)$, and commute with the outer $2$-sums. They therefore give $A_i\oplus_2B_j$ as a minor of $M$. Since $|E(H)|\geq3$, at least one element of $M$ is removed, so this minor is proper. Applying the same operations inside $A_j$ or $B_i$
gives the required minors retaining their roots, since the parallel pair transfers the root
from one end of $H$ to the other.
\end{proof}

For a rooted matroid $N$, let $v(N)\in\{0,1\}^{13}$ have entry one
in position $j$ exactly when $\pw(N\oplus_2T_j)\leq2$, and let
$w(N)$ be its number of ones. We also write $v(t)$ and $w(t)$
for the vector and its number of ones in row $t$ of
Table~\ref{tab:interface-codes}. The rooted matroid obtained from either component after deleting a decomposition edge is a proper minor of $M$, up to renaming its root. To obtain it, reduce the rooted matroid obtained from the other component to $U_{1,2}$ on its root and one element of the original ground set. Each of these two rooted matroids has at
least three elements, by the hypothesis on the sizes of the matroid labels, so the reduction
removes at least one element. Thus, in an excluded minor, the rooted
matroid obtained from either component has path-width at most two. It cannot belong to
the equivalence class corresponding to row $0$, since its sum with $U_{1,2}$ has path-width at most two whereas a rooted matroid in that class has no connected $2$-sum of path-width at most two. It cannot belong to the class corresponding to row $21$, whose members are equivalent to $U_{1,2}$, since its sum with the opposite side would then have path-width at most two.
The possible row indices for these equivalence classes therefore lie in
\[
 I=\{1,\ldots,20\}\setminus\{10,17\},\qquad 1\leq w(t)\leq12,
\]
using the indices of Table~\ref{tab:interface-codes}. For $a,b\in I$,
let $K(a,b)$ be one when a listed sequence for row $a$, or its reversal, admits an interleaving of width at most two with one for row $b$, or its reversal, and zero otherwise. Compute this value by \eqref{eq:interface-cell}. By Theorem~\ref{thm:finite-interface}, this is one exactly when the $2$-sum of rooted matroids in the corresponding classes has path-width at most two. This includes the graphic realization of the class corresponding to row $3$ in Proposition~\ref{prop:class-three}.

\begin{lemma}[A potential for successive cuts]
\label{lem:interface-potential}
For $1\leq r\leq12$, put
\[
 F(r)=13-r-\mathbf1_{r\leq7}-\mathbf1_{r\leq3},\qquad
 G(r)=r-\mathbf1_{r\geq5}-\mathbf1_{r\geq9},
\]
and $P(a,b)=\min\{F(w(a)),G(w(b))\}$. For $a,b,c,d\in I$, if
\[
 w(c)<w(a),\quad w(b)<w(d),\quad
 K(a,b)=K(c,d)=0,\quad K(a,d)=1,
\]
then $1\leq P(a,b)<P(c,d)\leq10$.
\end{lemma}

\begin{proof}
The values of $K$ obtained from Table~\ref{tab:interface-codes},
computed by \eqref{eq:interface-cell}, give
the following five finite facts for $x,y\in I$:
\begin{align}
 w(x)+w(y)\geq15&\ \Longrightarrow\ K(x,y)=1,\label{eq:interface-weight-sum}\\
 w(y)\geq7&\ \Longrightarrow\ K(11,y)=1,\label{eq:interface-type-eleven}\\
 w(y)\geq11&\ \Longrightarrow\ K(4,y)=1,\label{eq:interface-type-four}\\
 K(x,7)=1&\ \Longrightarrow\ w(x)\geq10,\label{eq:interface-type-seven}\\
 K(x,14)=1&\ \Longrightarrow\ w(x)\geq6.\label{eq:interface-type-fourteen}
\end{align}
The verifier described in Section~\ref{sec:certificates} reconstructs
these values from the sequences. The potential functions have values
\[
\begin{array}{c|rrrrrrrrrrrr}
r&1&2&3&4&5&6&7&8&9&10&11&12\\\hline
F(r)&10&9&8&8&7&6&5&5&4&3&2&1\\
G(r)&1&2&3&4&4&5&6&7&7&8&9&10.
\end{array}
\]
Both terms in the minimum are nondecreasing under the stated
transition. If the minimum stays unchanged, one term must remain
constant at that minimum. The displayed table shows two constant steps for $F$ and two for $G$, leaving four cases.

For $F$, the plateau from $w(a)=8$ to $w(c)=7$ has value five,
and $c=11$ is the only index in $I$ of weight seven. Equality requires
$G(w(b))\geq5$, hence $w(b)\geq6$.
By \eqref{eq:interface-weight-sum}, the equality $K(a,b)=0$ gives $w(b)\leq6$.
Then $w(d)\geq7$ contradicts $K(11,d)=0$ by
\eqref{eq:interface-type-eleven}.
The other plateau, from $w(a)=4$ to $w(c)=3$, has value eight
and forces $c=4$. The same argument gives $w(b)=10$ and
$w(d)\geq11$, contradicting \eqref{eq:interface-type-four}.

For $G$, the plateau from $w(b)=4$ to $w(d)=5$ has value four
and forces $d=7$. Equality requires $F(w(a))\geq4$, so
$w(a)\leq9$, contrary to $K(a,7)=1$ and
\eqref{eq:interface-type-seven}.
The remaining plateau, from $w(b)=8$ to $w(d)=9$, has value
seven and forces $d=14$. Equality gives $w(a)\leq5$, contrary
to \eqref{eq:interface-type-fourteen}. Hence the potential increases
strictly, and its range is given by the displayed table.
\end{proof}

\begin{theorem}[Decomposition diameter of excluded minors]
\label{thm:interface-diameter}
Let $T$ be the tree of a decomposition using $2$-sums that represents an excluded minor for path-width at most two. Suppose that every vertex is labelled by a connected matroid with at least three elements and that the distinguished elements for edges incident with the same vertex are distinct. Then $T$ has diameter at most ten edges. In particular, the tree in the canonical tree decomposition has at most eleven vertices on any path and can be rooted so that every vertex has distance at most five from the root.
\end{theorem}

\begin{proof}
Take a path with $m$ edges and rooted matroids $A_i,B_i$ obtained from its edge deletions as in
Lemma~\ref{lem:interface-splice}, with corresponding row indices $a_i,b_i$.
The original sums are rejected, so $K(a_i,b_i)=0$.
Shortening gives proper minors, so $K(a_i,b_j)=1$ for $i<j$.
Minor monotonicity gives $v(A_j)\leq v(A_i)$ and
$v(B_i)\leq v(B_j)$ coordinatewise. In each comparison, the two vectors must differ in at least one coordinate. Indeed, if $v(A_j)=v(A_i)$, Theorem~\ref{thm:finite-interface} makes $A_i\oplus_2B_j$ and $A_j\oplus_2B_j=M$ agree at width two, although the former is a proper minor and the latter is excluded. If $v(B_i)=v(B_j)$, comparing the same shortened sum with $A_i\oplus_2B_i=M$ gives the same contradiction. Thus
Lemma~\ref{lem:interface-potential} gives
\[
 1\leq P(a_1,b_1)<\cdots<P(a_m,b_m)\leq10,
\]
and $m\leq10$. Rooting at a centre of the tree gives depth at most five.
\end{proof}

Vertices labelled by $U_{1,2}$ must be excluded from this diameter bound, since inserting such a vertex would subdivide an edge without changing the represented matroid.
The theorem does not assert that diameter ten is attained.

\begin{corollary}\label{cor:interface-circuit-factors}
Every circuit or cocircuit labelling a vertex of the canonical tree decomposition of an excluded minor for path-width at most two has at most eleven elements, including the distinguished elements used for the incident $2$-sums.
\end{corollary}

\begin{proof}
Let a vertex be labelled by a circuit or cocircuit with $s$ elements and have degree at least two.
Replace this vertex by a path of $s-2$ vertices labelled by triangles, or by triads in the cocircuit case, placing two chosen distinguished elements at opposite ends. Extending
the path into the components beyond those two incident edges gives $(s-3)+2=s-1$ edges in a tree
whose vertices are labelled by connected matroids with at least three elements. Hence $s\leq11$ by
Theorem~\ref{thm:interface-diameter}.
When the vertex has degree one, the $s-1$ remaining elements form a series
or parallel class, so $s\leq3$ by Lemma~\ref{lem:loops-coloops}.
When the vertex has degree zero, the matroid itself would be a circuit
or cocircuit and have path-width at most one.
\end{proof}

\section{Complete binary classifications through rank five}\label{sec:low-rank}

We now restrict to binary matroids and classify the excluded minors
of rank at most five. No bound on the number of elements is imposed.
There are no excluded minors of rank at most two, since path-width
is at most rank.
The following criterion gives a short direct treatment of rank three.
Its binary sufficient direction occurs in
\cite[proof of Theorem~2.5]{LL16}.

\begin{lemma}[Cocircuit criterion]\label{thm:cocircuit-width}
A matroid $M$ of rank $r\geq2$ has $\pw(M)=r$ if and only if
every pair of cocircuits, allowing repetition, intersects in at
least two elements.
\end{lemma}

\begin{proof}
We have $\pw(M)<r$ exactly when two nonspanning sets and at most one
additional element cover $E(M)$. Indeed, in an ordering of width less than $r$, let $e$ complete the first spanning initial segment. The set before $e$ is nonspanning by this choice, and the cut immediately after $e$ has connectivity equal to the rank of its final segment, so that final segment is also nonspanning. Conversely, let $A,B$ be the two nonspanning sets in such a cover, and assign any overlap to one side. Listing $A$, then the possible additional element, and then $B$ gives an ordering in which every cut has a nonspanning side, so its width is less than $r$.
Extending the two sets to hyperplanes shows that such a cover exists
exactly when their complementary cocircuits intersect in at most
one element. Since path-width is at most rank, the result follows.
\end{proof}

\begin{theorem}[Rank three]\label{cor:rank-three-binary}
The rank-three binary excluded minors for path-width at most two are
exactly $F_7$ and $M(2K_4)$, where $2K_4$ has every edge doubled.
\end{theorem}

\begin{proof}
Let $M$ be such an excluded minor and put $N=\si(M)$.
Lemma~\ref{thm:parclass} gives $M=N[D]$. Every simple rank-three
binary matroid is a restriction of $F_7$. A restriction on at most
five points is covered by two nonspanning sets. To see this, choose
two omitted points $p,q$ and let $t$ be the third point on their
line. The other two lines through $t$ cover the remaining five
points. Replace each point of $N$ on these lines by its parallel class in $M$. The resulting two sets are nonspanning and cover $E(M)$, so Lemma~\ref{thm:cocircuit-width} rules out this case.
Thus $N$ is $F_7$ or $F_7\setminus p\cong M(K_4)$.

Distinct cocircuits of $F_7$ meet in two elements, so $\pw(F_7)=3$.
Every deletion is $M(K_4)$, whose vertex-star cocircuits have
singleton intersections and hence give path-width at most two. Every
contraction has rank two. Thus $F_7$ is an excluded minor, and no
proper parallel extension of it is an excluded minor.

For $N=M(K_4)$, the vertex stars at the ends of an edge meet in
that edge alone. The cocircuit criterion therefore forces every
edge to be doubled. Every bond of $K_4$ is either a vertex star or a cut separating two pairs of vertices. Two distinct vertex stars meet in their common edge, a vertex star meets such a cut in two edges, and two distinct cuts separating two pairs of vertices also meet in two edges. In the full double, every cocircuit
intersection consequently has size at least two. Deleting one copy
of an edge makes the intersection of its two endpoint stars a
singleton, giving path-width at most two. Every contraction has rank
two. Hence $M(2K_4)$ is the second excluded minor.
\end{proof}

\subsection{Graphic reduction}\label{sec:reduction}

To treat ranks four and five, we first reduce the problem to graphic
matroids and a short list of specified exceptions. Throughout,
$\mathcal{T}$ denotes the four excluded minors $F_7$, $F_7^{*}$, $M^{*}(K_5)$,
$M^{*}(K_{3,3})$ in Tutte's characterization. Every binary matroid with no minor in
$\mathcal{T}$ is graphic \cite{Tut59}.

\begin{lemma}[Six exceptional matroids]\label{lem:exceptional}
Each of the six exceptional matroids
\[
F_7,\quad F_7^{*},\quad M(K_5),\quad M^{*}(K_5),\quad
M(K_{3,3}),\quad M^{*}(K_{3,3})
\]
is an excluded minor for path-width at most two.
\end{lemma}

\begin{proof}
Theorem~\ref{cor:rank-three-binary} and duality give the two Fano
cases. For $R\in\{M(K_{3,3}),M(K_5)\}$, direct subset dynamic
programming gives
\[
 \pw(R)=3,\qquad \pw(R\setminus e)\leq2,\qquad \pw(R/e)\leq2
 \quad(e\in E(R)).
\]
The independent verifier checks every element from the stored binary
representations, as described in Section~\ref{sec:certificates}.
Duality gives the two cographic cases.
\end{proof}

\begin{lemma}[Reduction]\label{lem:reduction}
Every binary excluded minor for path-width at most two is one of the six
exceptional matroids in Lemma~\ref{lem:exceptional} or the cycle matroid $M(G)$
of a loopless $2$-connected multigraph $G$ on at least three vertices, with
$|E(G)| = |E(M)|$.
\end{lemma}

\begin{proof}
Let $M$ be an excluded minor not among these six. Each member of $\mathcal{T}$ has
path-width three by Lemma~\ref{lem:exceptional}, so a minor of $M$ lying in
$\mathcal{T}$ would be proper, contradicting that all proper minors of $M$ have
path-width at most two. Hence $M$ has no minor in $\mathcal{T}$ and is graphic by Tutte's theorem. Write $M=M(G)$, choosing $G$ connected and without isolated vertices. A loop of $G$ is a loop of $M$, excluded by
Lemma~\ref{lem:loops-coloops}. If $G$ has at most two vertices, then
$M\cong U_{1,m}$ and its ground set is a single parallel class. Its connectivity
function is bounded by one, so $\pw(M) \leq 1$ and $M$ is not an excluded minor.
So $G$ has at least three vertices. Finally $G$ is
$2$-connected. Otherwise $G$ has a cut vertex $v$, so $E(G)$ has a partition
$E_1 \cup E_2$ into nonempty sets such that the corresponding edge-induced
subgraphs intersect only in $v$. Every cycle of $G$ then lies wholly inside
$E_1$ or $E_2$, whence
$r_M(E_1) + r_M(E_2) = r_M(E)$ and $M = M|E_1 \oplus M|E_2$. By the direct-sum formula the path-width of $M$ is the maximum of the two parts,
so one part already has path-width three while being a proper minor of $M$,
contradicting minimality.
\end{proof}

Different graphs can have isomorphic cycle matroids, as described
by Whitney's $2$-isomorphism theorem~\cite{Whi33}. We therefore test
binary matroid isomorphism separately from graph isomorphism.
Fix a basis of the first matroid and consider every basis of the
second. For each basis $B$, write the representation as $[I_B\mid A]$.
Its \emph{fundamental graph} is bipartite, with one vertex for each
basis element and one for each element outside $B$. Join the vertices
corresponding to row $b$ and column $e$ exactly when $A_{b,e}=1$,
and give the two parts distinct colours. A colour-preserving graph
isomorphism identifies the two representations. Conversely, a
matroid isomorphism maps the fixed basis to one of the bases examined,
so this search detects every matroid isomorphism.

\subsection{Ranks four and five}

Let $M_2$ be the cycle matroid of the multigraph obtained from $K_5-uv$ by
doubling the three edges incident with $u$. This is the representative in Table~\ref{tab:multigraphs}, and
$|E(M_2)|=12$.

\begin{theorem}[Rank four]\label{thm:rank-four}
The rank-four binary excluded minors for path-width at most two are exactly
\[
F_7^*,\qquad M^*(K_{3,3}),\qquad M(K_5),\qquad M_2.
\]
In particular, every rank-four binary excluded minor has at most twelve
elements, and this bound is attained by $M_2$.
\end{theorem}

\begin{proof}
The four displayed matroids are excluded minors by
Lemma~\ref{lem:exceptional} and the direct verification of $M_2$.
For the converse, let $M$ be a rank-four binary excluded minor.
Lemma~\ref{lem:reduction} leaves the listed exceptional matroids or a
representation $M=M(G)$ by a loopless $2$-connected multigraph. Since
$r(M(G))=|V(G)|-1$, this graph has five vertices. Every parallel class has
at most two elements, so $G$ is obtained from a simple $2$-connected graph
$H$ on five vertices by doubling a set of edges. This reduction imposes
no bound on $|E(M)|$ in advance.

There are ten such simple graphs up to vertex isomorphism. The supplementary script generates all $2^{10}$ edge subsets on a fixed five-vertex set and tests $2$-connectivity directly. It then checks graph isomorphism by considering every vertex permutation, leaving ten representatives. Using the exact parallel-extension formula, it examines all $2{,}592$ doubling sets and finds exactly ten that are inclusion-minimal subject to producing a matroid of path-width greater than two. Seven yield a thirteen-element
matroid with a single-element contraction isomorphic to $M(2K_4)$, so
these cannot be excluded minors. Two give isomorphic copies of $M_2$ and
one gives $M(K_5)$. Thus the two displayed graphic matroids are exhaustive.
Section~\ref{sec:certificates} describes the independent checks of this search.
\end{proof}

Let $S_3$ be the graph with central triangle $012$ and three further
vertices $3,4,5$ adjacent respectively to $\{0,1\}$, $\{1,2\}$,
and $\{0,2\}$. Write $2S_3$ for the multigraph obtained by doubling
every edge of $S_3$.

The doubled $3$-sun is a known obstruction covered by \cite[Theorem~1]{KTY14}.
The next classification establishes its extremal uniqueness at rank five.

\begin{table}[htbp]
\centering
\small
\begin{tabular}{lrl}
\hline
matroid & vertices & edges \\
\hline
$M_2$ & 5 & $01^{2}\,02^{2}\,03^{2}\,12\,13\,14\,23\,24\,34$ \\
$M_3$ & 6 & $01^{2}\,02^{2}\,03^{2}\,12\,13\,14\,25\,35\,45$ \\
$M_4$ & 6 & $01^{2}\,02\,03\,14\,15\,23\,24^{2}\,35^{2}\,45$ \\
$M_5$ & 6 & $01^{2}\,02^{2}\,03\,04\,13\,14\,23\,24\,35\,45$ \\
$M_6$ & 6 & $01^{2}\,02\,03\,04\,12\,13\,15\,24\,25\,34\,35$ \\
$M_7$ & 6 & $01\,02\,03\,04\,12\,13\,15\,24\,25\,34\,35\,45$ \\
$L_3^{*}$ & 6 & $01^{2}\,02^{2}\,04\,05\,13^{2}\,14\,23^{2}\,24\,34\,45$ \\
$L_4^{*}$ & 6 & $01^{2}\,02^{2}\,04\,13^{2}\,14\,23^{2}\,24\,35^{2}\,45$ \\
$P_{12}$ & 6 & $01^{2}\,02\,03^{2}\,12\,14^{2}\,23\,24\,25\,35^{2}\,45^{2}$ \\
$P_{13}$ & 6 & $01^{2}\,02^{2}\,03^{2}\,12^{2}\,14\,25^{2}\,34^{2}\,45^{2}$ \\
\hline
\end{tabular}
\caption{The named graphic representatives used at ranks four and five.
The symbol $ij^2$ denotes two parallel edges. The superscript $*$ in
a matroid name denotes duality. The labels preserve those of the
verification data \cite{Kal26}. The matroids $M_2,\ldots,M_7$
already occur in Kashyap's partial list \cite[Fig.~4.3]{Kas08}.}
\label{tab:multigraphs}
\end{table}

\begin{theorem}[Rank five]\label{thm:rank-five}
The rank-five binary excluded minors for path-width at most two are exactly
\[
M(K_{3,3}),\quad M_3,M_4,M_5,M_6,M_7,\quad
L_3^*,L_4^*,\quad P_{12},P_{13},\quad M(2S_3),
\]
where the named representatives are defined in
Table~\ref{tab:multigraphs}.
The matroid $M(2S_3)$ has eighteen elements, and every other rank-five binary excluded minor has at most fifteen.
\end{theorem}

\begin{proof}
Among the six exceptional matroids, only $M(K_{3,3})$ has rank five,
and it is graphic. Thus every rank-five binary excluded minor has a representation $M=M(G)$, where $G$ is a loopless $2$-connected multigraph on six vertices. Since each parallel class has at most two elements, $G$ is obtained by doubling a set of edges of a simple $2$-connected graph on those vertices.

Generate all $2^{15}=32{,}768$ simple graphs on a fixed six-vertex set,
test connectivity after every vertex deletion, and identify graphs under
all vertex permutations, leaving $56$ representatives. We use no planarity filter.
The $56$ simple graphs have $124{,}608$ doubling subsets in total. For each simple graph,
Theorem~\ref{thm:exact-doubling} decides whether the parallel extension has path-width
at most two. If a proper doubling subset already gives a matroid of path-width greater than two, the parallel extension cannot be minor-minimal, since deleting a parallel copy at an edge outside that subset leaves the smaller parallel extension as a minor. Encode a doubling set by a binary integer whose one-bits mark its
edges, and process these integers in increasing order. Every proper
subset has a smaller integer, so a candidate discarded for this
reason has an earlier subset already known to produce a matroid of path-width greater
than two.

The $26{,}413$ remaining exact width decisions leave $113$
inclusion-minimal candidates. Deleting and contracting a representative
of each parallel class tests all distinct single-element minors.
Twelve multigraph isomorphism classes survive. Exact binary matroid
isomorphism identifies two of them with $L_4^*$, leaving the eleven
isomorphism classes in the statement. Ten are the nine rank-five matroids in Table~\ref{tab:multigraphs}
and $M(K_{3,3})$. The remaining isomorphism class is represented by $M(2S_3)$.

The listed sizes prove the extremal assertions.
Section~\ref{sec:certificates} describes the independent checks of this search.
\end{proof}

Table~\ref{tab:fixed-rank} summarizes the complete fixed-rank searches.

\begin{table}[htbp]
\centering
\small
\begin{tabular}{lrr}
\hline
quantity & rank four & rank five \\
\hline
simple $2$-connected graph isomorphism classes & 10 & 56 \\
doubling subsets covered & $2{,}592$ & $124{,}608$ \\
inclusion-minimal candidates & 10 & 113 \\
multigraph isomorphism classes retained & 2 & 12 \\
cycle-matroid isomorphism classes retained & 2 & 11 \\
binary excluded-minor isomorphism classes & 4 & 11 \\
maximum number of elements & 12 & 18 \\
\hline
\end{tabular}
\caption{Complete fixed-rank searches. The element maxima are conclusions
of the calculations, not input cutoffs. At rank four the two nongraphic
isomorphism classes are $F_7^*$ and $M^*(K_{3,3})$.}
\label{tab:fixed-rank}
\end{table}

There are therefore exactly $2+4+11=17$ isomorphism classes of binary excluded minors of rank at most five. The intersection with their dual list consists
of $F_7,F_7^*,M(K_{3,3}),M^*(K_{3,3})$, so their union with the
duals has thirty isomorphism classes. Every excluded minor outside this union
has rank and corank at least six.

\section{Computational verification}\label{sec:certificates}

Supplementary Material S1 contains the programs, input data,
transcripts and independent verification reports. All calculations
use exact integer arithmetic. The finite table in
Theorem~\ref{thm:finite-interface} and the fixed-rank classifications
require exhaustive calculations on the specified finite sets.
The general $2$-sum recurrence is proved above, and finite examples
check its implementation.

To reproduce these calculations, extract S1 into one directory and
follow the commands in \texttt{README\_S1.txt}. The Python programs
require Python 3.10 or later and only its standard library.
The rank-five search uses a C++17 compiler with GNU integer
bit-operation built-ins. In the archived program documentation,
\texttt{word} and \texttt{profile} refer to the sequences of triples
used here, and \texttt{language} refers to the set $\mathcal L$ of sequences admitting the specified interleavings. The program filenames and data labels are
retained so that the published calculations remain directly reproducible.

For each fixed-rank search, the transcript lists every simple graph and
minimal candidate. Independent graph generation recovers all ten simple graphs
at rank four and all fifty-six simple graphs at rank five. For every simple graph, subset ranks
computed from graph components agree with Gaussian elimination over $\F$. For rank four, the verifier checks the contractions that reject seven candidates. For rank five, it recomputes all $26{,}413$ unpruned width decisions by integer-valued dynamic programming and checks a specified proper minor of path-width greater than two for every rejected minimal candidate. Each positive rank-five candidate is matched to exactly one of the twelve multigraph representatives by comparing all vertex relabellings and giving an explicit vertex map. At both ranks, calculations
on the expanded ground sets check the widths and all single-element minors
of every passing representative, including $M(2S_3)$. Exact binary matroid
isomorphism then identifies the eleven rank-five classes. Together, these
checks cover the complete fixed-rank reductions, including rejected
candidates and positive representatives.

The verification of the finite table enumerates all $114$ sequences
and checks for each pair whether the grid has the required path with cut values at most two. It recomputes the sets $\mathcal L$ obtained from families closed under reversal and the restrictions forced by connected
deletion or contraction, leaving at most twenty sets.
Independent binary-rank calculations verify all twenty graphic or
cographic realizations. For the rooted matroid representing row $3$, a separate
component-count implementation checks every subset rank, the complete
family of sequences, and the path-widths of the series extension, the parallel extension and the $2$-sum with $T_{11}$ in
Proposition~\ref{prop:class-three}. The thirteen specified tests distinguish
all twenty possible sets. This finite table is part of the
proof of Theorem~\ref{thm:finite-interface}. The same grid-path calculation gives the fifteen distinguishing sets displayed in
Corollary~\ref{cor:minimum-tests}, which prove that thirteen fixed
tests are necessary. The program \texttt{check\_minimum\_tests.py}
checks these sets independently and writes
\texttt{minimum\_tests\_report.json}. Both files are included in S1.
The compression argument establishes the
applicability of the finite table to arbitrary matroids.

The additional program \texttt{check\_interface\_diameter.py}
reconstructs the values of $K$ from the listed sequences, verifies the
five finite facts used in Lemma~\ref{lem:interface-potential}, and
checks its potential on all $4{,}549$ permitted transitions among
$175$ ordered pairs for which $K=0$. The output
\texttt{interface\_diameter\_report.json} also gives maximum
sequence length ten for these formal transitions. These files are
included in Supplementary Material S1, archived in Zenodo
version~8~\cite{Kal26}. The minor operations needed for
Theorem~\ref{thm:interface-diameter} are proved above.

The program \texttt{canonical\_two\_sum.py} implements
Theorem~\ref{thm:general-two-sum-tree} for connected matroid labels supplied by rank oracles, i.e., procedures that return the rank of any specified subset. It retains full reduced families of sequences and the paths needed to reconstruct an ordering. After accepting at the root vertex, it expands the retained order of moves by feasible grid interleavings and checks every cut
of the resulting element ordering. Intermediate families use the
containing set of $208$ sequences and are not restricted to the $114$
sequences of actual orderings. Explicit resource limits report
exhaustion separately from a negative decision on whether path-width is at most two.

An independent checker tests $34$ matroids constructed by $2$-sums under $78$ choices of root vertex in their specified trees and checks all $73$ returned orderings directly against expanded rank tables. The cases include trees whose vertices are labelled by circuits, cocircuits and $3$-connected matroids, trees with vertices of degree greater than two or vertices labelled by $U_{1,2}$, and six examples involving nonbinary uniform matroids. For the binary cases, an independent representation of the complete
sum is constructed by taking vector-space quotients at the shared
distinguished elements. For the other cases, the ranks are computed
directly from the $2$-sum formula. Subset dynamic programming on the
complete ground set then determines the widths independently of the
tree recurrence. The examples have at most eighteen elements, and
these reference calculations examine $685{,}348$ subsets altogether. Separate programs check the
parallel-extension formula, the insertion formula for a
$3$-connected summand, and the endpoint counterexamples. The supplied commands and file hashes make the
scope of each calculation explicit. No polynomial running-time claim is made when the sizes of the matroid labels or the degrees of the tree vertices are unrestricted.

\section{Concluding remarks}\label{sec:outlook}

Whenever a connected rooted matroid occurs as a summand of a $2$-sum with another connected rooted matroid, it can be replaced by one of the twenty graphic or cographic representatives without changing whether the resulting path-width is at most two. The thirteen fixed matroids determine its equivalence class and hence which representative may be used. To combine the matroids obtained from several child subtrees, the recurrence retains full families of reduced sequences. For a $3$-connected summand on at least four elements, a consecutive ordering reduces replacement to five simpler conditions.

The fourteen-element representative of the equivalence class corresponding to row $3$ has series and parallel extensions of path-width two, but its $2$-sum with $T_{11}$ has path-width three. For $M(K_4)$ rooted at any edge, which represents row $14$, the series extension, the parallel extension and the $2$-sum with $T_{11}$ all have path-width at most two. Thus the series and parallel extension conditions do not distinguish these two equivalence classes. Corollary~\ref{cor:minimum-tests} determines the exact number of fixed connected rooted matroids needed for the complete classification, i.e., thirteen are both sufficient and necessary to distinguish all twenty equivalence classes.

For an excluded minor, the finite-test classification bounds the diameter of the tree in its canonical tree decomposition by ten edges. Every circuit or cocircuit labelling a vertex has at most eleven elements, including the distinguished elements used for the incident $2$-sums. These conclusions hold without a representability assumption. The bounds leave the number of vertices and the sizes of the $3$-connected matroids labelling them unrestricted. The realization of the class corresponding to row $3$ precludes improving these bounds merely by removing that class from the transition table. This realization does not establish that either bound is sharp. Smaller bounds may require further restrictions on the equivalence classes that can occur together along a decomposition path of an excluded minor.

The binary excluded-minor lists are complete through rank five and, by duality, through corank five. Further progress requires control of excluded minors with rank and corank both at least six. The criteria for the path-width of $2$-sums describe their possible tree decompositions, but the present results do not provide a global bound on the number of elements or a complete binary excluded-minor list.

\section*{Data and code availability}

All data, code and verification certificates supporting this article are
publicly available in the Zenodo v8 dataset~\cite{Kal26}
(Supplementary Material S1). The archive includes the complete fixed-rank
searches, reproduction instructions and file hashes.

\begingroup
\small
\setlength{\bibsep}{1pt}
\bibliographystyle{plainnat}
\bibliography{mybib}
\endgroup

\end{document}